\documentclass[a4paper,cleveref, autoref,numberwithinsect]{lipics-v2021}
\usepackage{graphicx}
\usepackage{color}
\usepackage{amsmath}
\usepackage{amssymb}
\usepackage{xspace}
\usepackage{epsfig}
\usepackage[dvipsnames]{xcolor}
\usepackage{soul}
\usepackage{xfrac}
\usepackage{enumitem}
\nolinenumbers
\newcommand {\mm}[1] {\ifmmode{#1}\else{\mbox{\(#1\)}}\fi}

\newcommand {\ceiling}[1] {{\left\lceil #1 \right\rceil}}

\newcommand{\ignore}[1]{}

\newsavebox{\smallProofsym}                 
\savebox{\smallProofsym}                               %
{
\begin{picture}(6,6)
\put(0,0){\framebox(6,6){}}
\put(0,2){\framebox(4,4){}}
\end{picture} 
}  

\makeatletter
\long\def\@makecaption#1#2{%
  \vskip\abovecaptionskip
  \sbox\@tempboxa{\small #1: #2}%
  \ifdim \wd\@tempboxa >\hsize
    \small #1: #2\par
  \else
    \global \@minipagefalse
    \hb@xt@\hsize{\hfil\box\@tempboxa\hfil}%
  \fi
  \vskip\belowcaptionskip}
\makeatother

\theoremstyle{plain}
\newtheorem*{namedthm}{\namedthmname}
\newcounter{namedthm}

\makeatletter
\newenvironment{named}[1]
  {\def\namedthmname{#1}%
   \refstepcounter{namedthm}%
   \namedthm\def\@currentlabel{#1}}
  {\endnamedthm}
\makeatother

\theoremstyle{plain}

\newtheorem*{supermaintheorem*}{Main Theorem}

\newtheorem*{supermaincorollary*}{Main Corollary}

\newcommand{\Rspace}        {\mm{{\mathbb R}}}
\newcommand{\Cech}[2]       {\mm{{\rm \v{C}ech}{({#1},{#2})}}}
\newcommand{\VRips}[2]      {\mm{{\rm Rips}{({#1},{#2})}}}
\newcommand{\Voronoi}[1]    {\mm{{\rm Vor}{({#1})}}}
\newcommand{\domain}[2]     {\mm{{\rm dom}{({#1},{#2})}}}
\newcommand{\Delaunay}[1]   {\mm{{\rm Del}{({#1})}}}
\newcommand{\Alf}[2]        {\mm{{\rm Alf}{({#1},{#2})}}}
\newcommand{\Alfa}[2]       {\mm{{\mathcal A}_{{#1},{#2}}}}
\newcommand{\RadiusF}       {\mm{{\rm Rad}}}
\newcommand{\ccz}[1]        {\mm{{z}_{#1}}}
\newcommand{\touch}[1]      {\mm{{\ell}{({#1})}}}
\newcommand{\short}[1]      {\mm{{j}{({#1})}}}
\newcommand{\Sphere}[1]     {\mm{S}^{#1}}

\newcommand{\Simplex}[1]    {\mm{\Sigma}_{#1}}

\newcommand{\radius}[3]     {\mm{r^{#1}_{{#2},{#3}}}}
\newcommand{\Circle}[1]     {\mm{{C}_{#1}}}
\newcommand{\radiusC}       {\mm{\gamma}}

\newcommand{\Betti}[2]      {\mm{\beta}_{#1}{#2}}

\newcommand{\conv}[1]       {\mm{{\rm conv\,}{#1}}}
\newcommand{\aff}[1]        {\mm{{\rm aff\,}{#1}}}
\newcommand{\dime}[1]       {\mm{\rm dim\,}{#1}}

\newcommand{\norm}[1]       {\mm{\|{#1}\|}}
\newcommand{\Edist}[2]      {\mm{\|{#1}-{#2}\|}}
\newcommand{\ee}            {\mm{\varepsilon}}

\newcommand{\Skip}[1]       {}

\title{Maximum Betti Numbers of \v{C}ech Complexes}
\titlerunning{Maximum Betti Numbers of \v{C}ech Complexes}

\author{Herbert Edelsbrunner}{ISTA (Institute of Science and Technology Austria), Kloster\-neu\-burg, Austria}{herbert.edelsbrunner@ist.ac.at}{https://orcid.org/0000-0002-9823-6833}{}

\author{J\'anos Pach}{R\'enyi Institute of Mathematics, Budapest, Hungary and ISTA (Institute of Science and Technology Austria), Kloster\-neu\-burg, Austria}{pach@cims.nyu.edu}{}{}

\authorrunning{Edelsbrunner, Pach}

\Copyright{Edelsbrunner, Pach}

\ccsdesc[100]{Theory of computation~Computational geometry}

\keywords{Discrete geometry, computational topology, \v{C}ech complexes, Delaunay mosaics, Alpha complexes, Betti numbers, extremal questions.}

\funding{\footnotesize
  The first author is supported by the European Research Council (ERC), grant no.\ 788183, and by the DFG Collaborative Research Center TRR 109, Austrian Science Fund (FWF), grant no.\ {I 02979-N35.}
  The second author is supported by the European Research Council (ERC), grant ``GeoScape'' and by the Hungarian Science Foundation (NKFIH), grant K-131529. Both authors are supported by the Wittgenstein Prize, Austrian Science Fund (FWF), grant no.\ Z 342-N31.
}

\begin{document}
\maketitle

\begin{abstract}
  The Upper Bound Theorem for convex polytopes implies that the $p$-th Betti number of the \v{C}ech complex of any set of $N$ points in $\Rspace^d$ and any radius satisfies $\Betti{p}{} = O(N^{m})$, with $m = \min \{ p+1, \ceiling{d/2} \}$.
  We construct sets in even and odd dimensions that prove this upper bound is asymptotically tight.
  For example, we describe a set of $N = 2(n+1)$ points in $\Rspace^3$ and two radii such that the first Betti number of the \v{C}ech complex at one radius is $(n+1)^2 - 1$, and the second Betti number of the \v{C}ech complex at the other radius is $n^2$.
\end{abstract}

\section{Introduction}
\label{sec:1}

Given a finite set of points $A$ in $\Rspace^d$ and a radius, their \emph{\v{C}ech complex} is the collection of subsets of the points whose balls have a nonempty common intersection.
This is an abstract simplicial complex isomorphic to the nerve of the balls, and by the Nerve Theorem~\cite{Bor48}, it has the same homotopy type as the union of the balls.
This property is the reason for the popularity of the \v{C}ech complex in topological data analysis; see e.g.~\cite{Car09,EdHa10}.
Of particular interest are the \emph{Betti numbers}, which may be interpreted as the numbers of holes of different dimensions.
These are intrinsic properties, but for a space embedded in $\Rspace^d$, they describe the connectivity of the space as well as that of its complement.
Most notably, the (reduced) zero-th Betti number, $\Betti{0}{}$, is one less than the number of \emph{connected components}, and the last possibly non-zero Betti number, $\Betti{d-1}{}$, is the number of \emph{voids} (bounded components of the complement).
Spaces that have the same homotopy type---such as a union of balls and the corresponding \v{C}ech complex---have identical Betti numbers.
While the \v{C}ech complex is not necessarily embedded in $\Rspace^d$, the corresponding union of balls is, which implies that also the \v{C}ech complex has no non-zero Betti numbers beyond dimension $d-1$.
To gain insight into the statistical behavior of the Betti numbers of \v{C}ech complexes, it is useful to understand how large the numbers can get, and this is the question we study in this paper.

\smallskip
The question of maximum Betti numbers lies at the crossroads of computational topology and discrete geometry.
Originally inspired by problems in the theory of polytopes ~\cite{McM71,Zie95}, optimization~\cite{Pad99},  robotics, motion planning~\cite{ScSh88}, and molecular modeling~\cite{Mez90}, many interesting and surprisingly difficult questions were asked about the complexity of the union of $n$ geometric objects, as $n$ tends to infinity. 
For a survey, consult~\cite{APS08}.
Particular attention was given to estimating the number of voids among $N$ simply shaped bodies, e.g., for the translates of a fixed convex body in $\Rspace^d$.
In the plane, the answer is typically linear in $N$ (for instance, for disks or other fat objects), but for $d=3$, the situation is more delicate. 
The maximum number of voids among $N$ translates of a convex polytope with a constant number of faces is $\Theta(N^2)$, but this number reduces to linear for the cube and other simple shapes~\cite{ArSh97}. 
It was conjectured for a long time that similar bounds hold for the translates of a convex shape that is not necessarily a polytope.
However, this turned out to be false: Aronov, Cheung, Dobbins and Goaoc~\cite{ACDG17} constructed a convex body in $\Rspace^3$ for which the number of voids is $\Omega(N^3)$. 
This is the largest possible order of magnitude for any arrangement of convex bodies that are not necessarily translates of a fixed one~\cite{Kov88}. 
It is an outstanding open problem whether there exists a \emph{centrally symmetric} convex body with this property. 

\smallskip
For the special case where the convex body is the \emph{unit ball} in $\Rspace^3$, the maximum number of voids in a union of $N$ translates is $O(N^2)$.
This can be easily derived from the Upper Bound Theorem for $4$-dimensional convex polytopes.
It has been open for a long time whether this bound can be attained.
Our main theorem answers this question in the affirmative, in a more general sense.
\begin{named}{Main Theorem}
  \label{thm:main}
  For every $d \geq 1$, $0 \leq p \leq d-1$, and $N \geq 1$, there is a set of $N$ points in $\Rspace^d$ and a radius such that the $p$-th Betti number of the \v{C}ech complex of the points and the radius is $\Betti{p}{} = \Theta (N^m)$, with $m = \min \{ p+1, \ceiling{d/2} \}$.
\end{named}
For $d=3$, the maximum second Betti number is $\Betti{2}{} = \Theta (N^2)$ in $\Rspace^3$, which is equivalent to the maximum number of voids being $\Theta (N^2)$.
In addition to the \v{C}ech complex, the proof of the Main Theorem makes use of three complexes defined for a set of $N$ points, $A \subseteq \Rspace^d$, in which the third also depends on a radius $r \geq 0$:
\begin{itemize}
  \item the \emph{Voronoi domain} of a point $a \in A$, denoted $\domain{a}{A}$, contains all points $x \in \Rspace^d$ that are at least as close to $a$ as to any other point in $A$, and the \emph{Voronoi tessellation} of $A$, denoted $\Voronoi{A}$, is the collection of $\domain{a}{A}$ with $a \in A$ \cite{Vor070809};
  \item the \emph{Delaunay mosaic} of $A$, denoted $\Delaunay{A}$, contains the convex hull of $\Sigma \subseteq A$ if the common intersection of the $\domain{a}{A}$, with $a \in \Sigma$, is non-empty, and no other Voronoi domain contains this common intersection \cite{Del34};
  \item the \emph{Alpha complex} of $A$ and $r$, denoted $\Alf{A}{r}$, is the subcomplex of the Delaunay mosaic that contains the convex hull of $\Sigma$ if the common intersection of the $\domain{a}{A}$, with $a \in \Sigma$, contain a point at distance at most $r$ from the points in $\Sigma$ \cite{EKS83,EdMu94}.
\end{itemize}
The Delaunay mosaic is also known as the \emph{dual} of the Voronoi tessellation, or the \emph{Delaunay triangulation} of $A$.
Note that $\Alf{A}{r} \subseteq \Alf{A}{R}$ whenever $r \leq R$, and that for sufficiently large radius, the Alpha complex is the Delaunay mosaic.
Similar to the \v{C}ech complex, the Alpha complex has the same homotopy type as the union of balls with radius $r$ centered at the points in $A$, and thus the same Betti numbers.
It is instructive to increase $r$ from $0$ to $\infty$ and to consider the \emph{filtration} or nested sequence of Alpha complexes.
The difference between an Alpha complex, $K$, and the next Alpha complex in the filtration, $L$, consists of one or more cells.
If it is a single cell of dimension $p$, then either $\Betti{p}{} (L) = \Betti{p}{} (L) + 1$ or $\Betti{p-1}{} (L) = \Betti{p-1}{} (L) - 1$, and all other Betti numbers are the same.
In the first case, we say the cell gives \emph{birth} to a $p$-cycle, while in the second case, it gives \emph{death} to a $(p-1)$-cycle, and in both cases we say it is \emph{critical}.
If there are two or more cells in the difference, this may be a generic event or accidental due to non-generic position of the points.
In the simplest generic case, we simultaneously add two cells (one a face of the other), and the addition is an anti-collapse, which does not affect the homotopy type of the complex.
More elaborate anti-collapses, such as the simultaneous addition of an edge, two triangles, and a tetrahedron, can arise generically.
The cells in an interval of size $2$ or larger cancel each other's effect on the homotopy type, so we say these cells are \emph{non-critical}.
We refer to \cite{BaEd17} for more details.

\smallskip
With these notions, it is not difficult to prove the upper bounds in the Main Theorem.
As mentioned above, the \v{C}ech and alpha complexes for radius $r$ have the same Betti numbers.
Since a $p$-cycle is given birth to by a $p$-cell in the filtration of Alpha complexes, and every $p$-cell gives birth to at most one $p$-cycle, the number of $p$-cells is an upper bound on the number of $p$-cycles, which are counted by the $p$-th Betti number.
The number of $p$-cells in the Alpha complex is at most that number in the Delaunay mosaic, which is, by the Upper Bound Theorem for convex polytopes \cite{McM71,Zie95}, at most $O(N^m)$, with $m = \min \{p+1, \ceiling{d/2}\}$.

\smallskip
By comparison, to come up with constructions that prove matching lower bounds is delicate and the main contribution of this paper.
Our constructions are multipartite and inspired by Lenz' constructions related to Erd\H os's celebrated question on repeated distances~\cite{Erd60}: what is the largest number of point pairs in an $N$-element set in $\Rspace^d$ that are at distance $1$ apart?
Lenz noticed that in $4$ (and higher) dimensions, this maximum is $\Theta (N^2)$. To see this, take two circles of radius $\sfrac{\sqrt{2}}{2}$ centered at the origin, lying in two orthogonal planes, and place $\lceil N/2\rceil$ and $\lfloor N/2\rfloor$ points on them.
By Pythagoras' theorem, any two points on different circles are at distance $1$ apart, so the number of unit distances is roughly $N^2/4$, which is nearly optimal. 
For $d=2$ and $3$, we are far from knowing asymptotically tight bounds.
The current best constructions give $\Omega (N^{1+c/\log \log N})$ unit distance pairs in the plane~\cite[page 191]{BMP05} and $\Omega(N^{4/3}\log\log N)$ in $\Rspace^3$, while the corresponding upper bounds are $O(N^{4/3})$ and $O(N^{3/2})$; see~\cite{SST84} and~\cite{KMSS12,Zah13}.
Even the following, potentially simpler, bipartite analogue of the repeated distance question is open in $\Rspace^3$: given $N$ red points and $N$ blue points in $\Rspace^3$, such that the minimum distance between a red and a blue point is $1$, what is the largest number of red-blue point pairs that determine a unit distance?
The best known upper bound, due to Edelsbrunner and Sharir~\cite{EdSh91} is $O(N^{4/3})$, but we have no superlinear lower bound. 
This last question is closely related to the subject of our present paper.

It is not difficult to see that the upper bounds in the Main Theorem also hold for the Betti numbers of the union of $N$ \emph{not necessarily congruent} balls in $\Rspace^d$. This requires the use of weighted versions of the Voronoi tessellation and the Upper Bound Theorem. 
In the lower bound constructions, much of the difficulty stems from the fact that we insist on using congruent balls. This suggests the analogy to the problem of repeated distances.

\medskip \noindent \textbf{Outline.}
Section~\ref{sec:2} proves the Main Theorem for sets in \emph{even} dimensions. 
Starting with Lenz' constructions, we partition the Delaunay mosaic into finitely many groups of \emph{congruent} simplices.
We compute the radii of their circumspheres and obtain the Betti numbers by straightforward counting.
In Section~\ref{sec:3}, we establish the Main Theorem for sets in three dimensions.
The situation is more delicate now, because the simplices of the Delaunay mosaic no longer fall into a small number of distinct congruence classes.
Nevertheless, they can be divided into groups of nearly congruent simplices, which will be sufficient to carry out the counting argument.
In Section~\ref{sec:4}, we extend the result to any \emph{odd} dimension. 
Again we require a detailed analysis of the shapes and sizes of the simplices, which now proceeds by induction on the dimension.
Section~\ref{sec:5} contains concluding remarks and open questions.

\section{Even Dimensions}
\label{sec:2}

In this section, we give an answer to the maximum Betti number question for \v{C}ech complexes in even dimensions.
To state the result, let $n_k$ be the minimum integer such that the edges of a regular $n_k$-gon inscribed in a circle of radius $1/\sqrt2$ are strictly shorter than $\sqrt{2/k}$.
For example, if $k=2$, we have $n_2=5$, as the side length of an inscribed square is equal to $1$.
\begin{theorem}[Maximum Betti Numbers in $\Rspace^{2k}$]
  \label{thm:maximum_Betti_numbers_in_even_dimensions}
  For every $2k \geq 2$ and $n \geq n_k$, there exist a set $A$ of $N = kn$ points in $\Rspace^{2k}$ and radii $\rho_0 < \rho_1 < \ldots < \rho_{2k-2}$ such that
  \begin{align}
    \Betti{p}{(\Cech{A}{\rho_p})} &= \tbinom{k}{p+1} \cdot n^{p+1} \pm O(1), \mbox{\rm ~~~for~} 0 \leq p \leq k-1 ; 
      \label{eqn:even1} \\
    \Betti{p}{(\Cech{A}{\rho_p})} &= \tbinom{k-1}{p+1-k} \cdot n^k \pm O(1), \mbox{\rm ~~~for~} k \leq p \leq 2k-2.
      \label{eqn:even2}
  \end{align}
  For $p = 2k-1$, there exist $N = k(n+1) + 2$ points in $\Rspace^{2k}$ and a radius such that the $p$-th Betti number of the \v{C}ech complex is $n^k \pm O(n^{k-1})$.
\end{theorem}
The reason for the condition $n \geq n_k$ will become clear in the proof of Lemma~\ref{lem:ordering_of_radii_in_even_dimensions}, which establishes a particular ordering of the circumradii of the cells in the Delaunay mosaic.
The proof of the cases $0 \leq p \leq 2k-2$ is not difficult using elementary computations, the results of which will be instrumental for establishing the more challenging odd-dimensional statements in Sections~\ref{sec:3} and \ref{sec:4}.
The proof consists of four steps presented in four subsections: the construction of the point set in Section~\ref{sec:2.1}, the geometric analysis of the simplices in the Delaunay mosaic in Section~\ref{sec:2.2}, the ordering of the circumradii in Section~\ref{sec:2.3}, and the final counting in Section~\ref{sec:2.4}.
The proof of the case $p = 2k-1$ in $\Rspace^{2k}$ readily follows the case $p = 2k-2$ in $\Rspace^{2k-1}$, as we will describe in Section~\ref{sec:4.5}.

\subsection{Construction}
\label{sec:2.1}

Let $d = 2k$.
We construct a set $A = A_{2k} (n)$ of $N = kn$ points in $\Rspace^d$ using $k$ concentric circles in mutually orthogonal coordinate planes:  for $0 \leq \ell \leq k-1$, the circle $\Circle{\ell}$ with center at the origin, $0 \in \Rspace^{d}$, is defined by $x_{2 \ell + 1}^2 + x_{2 \ell+2}^2 = \frac{1}{2}$ and $x_i = 0$ for all $i \neq 2\ell+1, 2\ell+2$.
On each of the $k$ circles, we choose $n \geq 3$ points that form a regular $n$-gon.
The length of the edges of these $n$-gons will be denoted by $2s$.
Obviously, we have $s = \frac{\sqrt{2}}{2} \sin \frac{\pi}{n}$.
Assuming $k \geq 2$, the condition $n \geq n_k$ implies that the Euclidean distance between consecutive points along the same circle is less than $1$, and by Pythagoras' theorem, the distance between any two points on different circles is $1$.
It follows that for $r = \frac{1}{2}$, neighboring balls centered on the same circle overlap, while the balls centered on different circles only touch.
Correspondingly, the first Betti number of the \v{C}ech complex for a radius slightly less than $\frac{1}{2}$ is $\Betti{1}{} = k$.
To get the first Betti number for $r = \frac{1}{2}$, we add all edges of length $1$, of which $k-1$ connect the $k$ circles into a single connected component, while the others increase the first Betti number to $\Betti{1}{} = k + \binom{k}{2} n^2 - (k-1) = \binom{k}{2} n^2 + 1$.

\smallskip
To generalize the analysis beyond the first Betti number, we consider the Delaunay mosaic and two radii defined for each of its cells.
The \emph{circumsphere} of a $p$-cell is the unique $(p-1)$-sphere that passes through its vertices, and we call its center and radius the \emph{circumcenter} and the \emph{circumradius} of the cell.
To define the second radius, we call a $(d-1)$-sphere \emph{empty} if all points of $A$ lie on or outside the sphere.
The \emph{radius function} on the Delaunay mosaic, $\RadiusF \colon \Delaunay{A} \to \Rspace$, maps each cell to the radius of the smallest empty $(d-1)$-sphere that passes through the vertices of the cell.
By construction, each Alpha complex is a sublevel set of this function: $\Alf{A}{r} = \RadiusF^{-1} [0,r]$.
The two radii of a cell may be different, but they agree for the critical cells as defined in terms of their topological effect in the introduction.
It will be convenient to work with the corresponding geometric characterization of criticality:
\begin{definition}[Critical Cell]
  \label{def:critical_cell}
  A \emph{critical cell} of $\RadiusF \colon \Delaunay{A} \to \Rspace$ is a cell $\Sigma \in \Delaunay{A}$ that (1) contains the circumcenter in its interior, and (2) the $(d-1)$-sphere centered at the circumcenter that passes through the vertices of $\Sigma$ is empty and the vertices of $\Sigma$ are the only points of $A$ on this sphere.
\end{definition}
There are two conditions for a cell to be critical for a reason.
The first guarantees that its topological effect is not canceled by one of its faces, and the second guarantees that it does not cancel the topological effect of one of the cells it is a face of.
As proved in \cite{BaEd17}, the radius function of a generic set, $A \subseteq \Rspace^d$, is \emph{generalized discrete Morse}; see Forman~\cite{For98} for background on discrete Morse functions.
This means that each level set of $\RadiusF$ is a union of disjoint combinatorial intervals, and a simplex is critical iff it is the only simplex in its interval.
Our set $A$ is not generic because the $(d-1)$-sphere with center $0 \in \Rspace^{2k}$ and radius $\sfrac{\sqrt{2}}{2}$ passes through all its points.
Indeed, $\Delaunay{A}$ is really a $2k$-dimensional convex polytope, namely the convex hull of $A$ and all its faces.
Nevertheless, the distinction between critical and non-critical cells is still meaningful, and all cells in the Delaunay mosaic of our construction will be seen to be critical.

\smallskip
The value of the $2k$-polytope under the radius function is $\sfrac{\sqrt{2}}{2}$, while the values of its proper faces are strictly smaller than $\sfrac{\sqrt{2}}{2}$.
Let $\Simplex{\ell,j}$ be such a face, in which $\ell+1$ is the number of circles that contain one or two of its vertices, and $j+1$ is the number of circles that contain two.
Specifically, $\Simplex{\ell,j}$ has $j+1$ disjoint \emph{short} edges of length $2s$, while the remaining \emph{long} edges all have unit length.
Indeed, the geometry of the simplex is determined by $\ell$ and $j$ and does not depend on the circles from which we pick the vertices or where along these circles we pick them, as long as two vertices from the same circle are consecutive along this circle.
For example, $\Simplex{1,-1}$, $\Simplex{1,0}$, and $\Simplex{1,1}$ are the unit length edge, the isosceles triangle with one short and two long edges, and the tetrahedron with two disjoint short and four long edges, respectively.
We call the $\Simplex{\ell,j}$ \emph{ideal simplices}.
In even dimensions they are \emph{precisely} the simplices in the Delaunay mosaic of our construction.
However, in odd dimensions, the cells in the Delaunay mosaic only \emph{converge} to the ideal simplices.
This will be explained in detail in Sections~\ref{sec:3} and \ref{sec:4}.

\subsection{Circumradii of Ideal Simplices}
\label{sec:2.2}

In this section, we compute the sizes of some ideal simplices, beginning in four dimensions.
The \emph{ideal $2$-simplex} or \emph{triangle}, denoted $\Simplex{1,0}$, is the isosceles triangle with one short and two long edges.
We write $h(s)$ for the \emph{height} of $\Simplex{1,0}$ (the distance between the midpoint of the short edge and the opposite vertex), and $r(s)$ for the circumradius.
There is a unique way to glue four such triangles to form the boundary of a tetrahedron: the two short edges are disjoint and their endpoints are connected by four long edges.
This is the \emph{ideal $3$-simplex} or \emph{tetrahedron}, denoted $\Simplex{1,1}$.
We write $H(s)$ for its \emph{height} (the distance between the midpoints of the two short edges), and $R(s)$ for its circumradius.
\begin{lemma}[Ideal Triangle and Tetrahedron]
  \label{lem:ideal_triangle_and_tetrahedron}
  The squared heights and circumradii of the ideal triangle and the ideal tetrahedron in $\Rspace^4$ satisfy
  \begin{align}
    h^2(s)  &=  1 - s^2 , ~~~~~~~~ 4r^2(s) = \frac{1}{1-s^2} , 
      \label{eqn:idealtriangle} \\
    H^2(s)  &=  1 - 2s^2,   ~~~~~~ 4R^2(s) = 1 + 2s^2 .
      \label{eqn:idealtetrahedron}
  \end{align}
\end{lemma}
\begin{proof}
  By Pythagoras' theorem, the squared height of the ideal triangle is $h^2 = 1-s^2$.
  If we glue the two halves of a scaled copy of the ideal triangle to the two halves of the short edge, we get a quadrangle inscribed in the circumcircle of the triangle.
  One of its diagonals passes through the center, and its squared length satisfies $4r^2  =  1 + (s/h)^2  =  1 + \frac{s^2}{1-s^2}$.

   By Pythagoras' theorem, the squared height of the ideal tetrahedron is $H^2 = h^2 - s^2 = 1 - 2 s^2$.
   Hence, the squared diameter of the circumsphere is $4R^2 = H^2 + (2s)^2 = 1+2s^2$.
\end{proof}
To generalize the analysis beyond the ideal simplices in four dimensions, we write $\radius{}{\ell}{j} (s)$ for the circumradius of $\Simplex{\ell,j}$, so $\radius{}{1}{-1} (s) = \frac{1}{2}$, $\radius{}{1}{0} (s) = r(s)$, and $\radius{}{1}{1} (s) = R(s)$.
For two kinds of ideal simplices, the circumradii are particularly easy to compute, namely for the $\Simplex{\ell,-1}$ and the $\Simplex{\ell,\ell}$, and we will see that knowing their circumradii will be sufficient for our purposes.
\begin{lemma}[Further Ideal Simplices]
  \label{lem:further_ideal_simplices}
  For $\ell \geq 0$, the squared circumradii of $\Simplex{\ell,-1}$ and $\Simplex{\ell,\ell}$ satisfy
  $\radius{2}{\ell}{-1} (s) = {\ell}/{(2\ell+2)}$ and $\radius{2}{\ell}{\ell} (s) = {(\ell+2s^2)}/{(2\ell+2)}$.
\end{lemma}
\begin{proof}
  Consider the standard $\ell$-simplex, which is the convex hull of the endpoints of the $\ell+1$ unit coordinate vectors in $\Rspace^{\ell+1}$.
  Its squared circumradius is the squared distance between the barycenter and any one of the vertices, which is easy to compute.
  By comparison, the squared circumradius of the regular $\ell$-simplex with unit length edges is half that of the standard $\ell$-simplex:
  \begin{align}
    R_\ell^2  &=  \frac{1}{2} \left[ \frac{\ell^2}{(\ell+1)^2} + \frac{1}{(\ell+1)^2} + \ldots + \frac{1}{(\ell+1)^2} \right]
              =  \frac{\ell}{2(\ell+1)} , 
      \label{eqn:regular_radius}
  \end{align}
  Since $\radius{2}{\ell}{-1} (s) = R_\ell^2$, this proves the first equation in the lemma.
  Note that the convex hull of the midpoints of the $\ell+1$ short edges of $\Simplex{\ell,\ell}$ is a regular $\ell$-simplex with edges of squared length $H^2(s) = 1-2s^2$.
  The short edges are orthogonal to this $\ell$-simplex, which implies
  \begin{align}
    \radius{2}{\ell}{\ell}  &=  H^2(s) \cdot R_\ell^2 + s^2
      =  R_\ell^2 + (1-2R_\ell^2) s^2
      =  \frac{\ell+2s^2}{2\ell+2} ,
  \end{align}
  which proves the second equation in the lemma.
\end{proof}

\subsection{Ordering the Radii}
\label{sec:2.3}

In this subsection, we show that the radii of the circumspheres of the ideal simplices increase with increasing $\ell$ and $j$:
\begin{lemma}[Ordering of Radii in $\Rspace^{2k}$]
    \label{lem:ordering_of_radii_in_even_dimensions}
    Let $0 < s < 1 / \sqrt{2k}$.
    Then the ideal simplices satisfy $\radius{}{\ell}{\ell} (s) < \radius{}{\ell+1}{-1} (s)$ for $0 \leq \ell \leq k-2$, and $\radius{}{\ell}{j} (s) < \radius{}{\ell}{j+1} (s)$ for $-1 \leq j < \ell \leq k-1$.
\end{lemma}
\begin{proof}
  To prove the first inequality, we use Lemma~\ref{lem:further_ideal_simplices} to compute the difference between the two squared radii:
  \begin{align}
    \radius{2}{\ell+1}{-1} (s) - \radius{2}{\ell}{\ell} (s)
      &=  \frac{\ell+1}{2(\ell+2)} - \frac{\ell+2s^2}{2(\ell+1)}
       =  \frac{1 - 2s^2(\ell+2)}{2 (\ell+2) (\ell+1)} .
  \end{align}
  Hence, $\radius{2}{\ell}{\ell} (s) < \radius{2}{\ell+1}{-1} (s)$ iff $s^2 < 1 / (2\ell+4)$.
  We need this inequality for $0 \leq \ell \leq k-2$, so $s^2 < 1 / (2k)$ is sufficient, but this is guaranteed by the assumption.

  \smallskip
  We prove the second inequality geometrically, without explicit computation of the radii.
  Fix an ideal simplex, $\Simplex{\ell,j}$, and let $\Sphere{d-1}$ be the $(d-1)$-sphere whose center and radius are the circumcenter and circumradius of $\Simplex{\ell,j}$.
  Assume w.l.o.g.\ that the circles $C_0$ to $C_j$ contain two vertices of $\Simplex{\ell,j}$ each, and the circles $C_{j+1}$ to $C_\ell$ contain one vertex of $\Simplex{\ell,j}$ each.
  For $0 \leq i \leq k-1$, write $P_i$ for the $2$-plane that contains $C_i$ and $x_i$ for the projection of the center of $\Sphere{d-1}$ onto $P_i$.
  Note that $\norm{x_i}^2$ is the squared distance to the origin, and for $0 \leq i \leq \ell$ write $r_i^2$ for the squared distance between $x_i$ and the one or two vertices of $\Simplex{\ell,j}$ in $P_i$.
  Fixing $i$ between $0$ and $\ell$, the squared radius of $\Sphere{d-1}$ is $r_i^2$ plus the squared distance of the center of $\Sphere{d-1}$ from $P_i$, which is the sum of the squared norms other than $\norm{x_i}^2$.
  Taking the sum for $0 \leq i \leq \ell$ and dividing by $\ell+1$, we get
  \begin{align}
    \radius{2}{\ell}{j} (s)  &=  \frac{1}{\ell+1} 
    \left[ \sum\nolimits_{i=0}^{\ell} r_i^2 + \ell \cdot \sum\nolimits_{i=0}^{\ell} \norm{x_i}^2 + (\ell+1) \cdot \sum\nolimits_{i={\ell+1}}^{k-1} \norm{x_i}^2 \right] .
      \label{eqn:minradius}
  \end{align}
  By construction, $\radius{2}{\ell}{j} (s)$ is the minimum squared radius of any $(d-1)$-sphere that passes through the vertices of $\Simplex{\ell,j}$.
  Hence, also the right-hand side of \eqref{eqn:minradius} is a minimum, but since the $2$-planes are pairwise orthogonal, we can minimize in each $2$-plane independently of the other.
  For $\ell+1 \leq i \leq k-1$, this implies $\norm{x_i}^2 = 0$, so we can drop the last sum in \eqref{eqn:minradius}.
  For $j+1 \leq i \leq \ell$, $x_i$ lies on the line passing through the one vertex in $P_i$ and the origin.
  This implies that $\Sphere{d-1}$ touches $C_i$ at this vertex, and all other points of the circle lie strictly outside $\Sphere{d-1}$.
  For $0 \leq i \leq j$, $x_i$ lies on the bisector line of the two vertices, which passes through the origin.
  The contribution to \eqref{eqn:minradius} for an index between $0$ and $j$ is thus strictly larger than for an index between $j+1$ and $\ell$.
  This finally implies $\radius{2}{\ell}{j} (s) < \radius{2}{\ell}{j+1} (s)$ and completes the proof of the second inequality.
\end{proof}

Recall that $2s$ is the edge length of a regular $n$-gon inscribed in a circle of radius $1/\sqrt{2}$. By the definition of $n_k$, the condition $s < 1 / \sqrt{2k}$ in the lemma holds,  whenever $n \geq n_k$.

\smallskip
For the counting argument in the next subsection, we need the ordering of the radii as defined by the radius function, but it is now easy to see that they are the same as the circumradii, so Lemma~\ref{lem:ordering_of_radii_in_even_dimensions} applies.
Indeed, $\RadiusF (\Simplex{\ell,j}) = \radius{}{\ell}{j} (s)$ if $\Simplex{\ell,j}$ is a critical simplex of $\RadiusF$.
To realize that it is, we note that the circumcenter of $\Simplex{\ell,j}$ lies in its interior because of symmetry.
To see that also the second condition for criticality in Definition~\ref{def:critical_cell} is satisfied, we recall that $\Sphere{d-1}$ is the $(d-1)$-sphere whose center and radius are the circumcenter and circumradius of $\Simplex{\ell,j}$.
By the argument in the proof of Lemma~\ref{lem:ordering_of_radii_in_even_dimensions}, $\Sphere{d-1}$ is empty, and all points of $A$ other than the vertices of $\Simplex{\ell,j}$ lie strictly outside this sphere.

\subsection{Counting the Cycles}
\label{sec:2.4}

To compute the Betti numbers, we make essential use of the structure of the Delaunay mosaic of $A$, which consists of as many groups of congruent ideal simplices as there are different values of the radius function.
For each $0 \leq \ell \leq k-1$, we have $\ell+2$ groups of simplices that touch exactly $\ell+1$ of the $k$ circles.
In addition, we have a single $2k$-cell, $\conv{A}$, with radius $\sfrac{\sqrt{2}}{2}$, which gives $1 + 2 + \ldots + (k+1) = \binom{k+2}{2}$ groups.
We write $\Alfa{\ell}{j} = \RadiusF^{-1} [0, \radius{}{\ell}{j}]$ for the Alpha complex that consists of all simplices with circumradii up to $\radius{}{\ell}{j} = \radius{}{\ell}{j} (s)$.
We prove Theorem~\ref{thm:maximum_Betti_numbers_in_even_dimensions} in two steps, first the relations \eqref{eqn:even1} for $0 \leq p \leq k-1$ and second the relations \eqref{eqn:even2} for $k \leq p \leq 2k-2$.
The case $p = 2k-1$ will be settled later, in Section~\ref{sec:4.5}.
To begin, we study the Alpha complexes whose simplices touch at most $\ell+1$ of the $k$ circles.
\begin{lemma}[Constant Homology in $\Rspace^{2k}$]
  \label{lem:constant_homology_in_even_dimensions}
  Let $k$ be a constant, $A = A_{2k} (n) \subseteq \Rspace^{2k}$, and $0 \leq \ell \leq k-1$.
  Then $\Betti{p}{(\Alfa{\ell}{\ell})} = O(1)$ for every $0 \leq p \leq 2k-1$.
\end{lemma}
\begin{proof}
  Fix $\ell$ and a subset of $\ell+1$ circles.
  The full subcomplex of $\Alfa{\ell}{\ell}$ defined by the points of $A$ on these $\ell+1$ circles consists of all cells in $\Delaunay{A}$ whose vertices lie on these and not any of the other circles.
  Its homotopy type is that of the join of $\ell+1$ circles or, equivalently, that of the $(2\ell+1)$-sphere; see \cite[pages 9 and 19]{Hat02}.
  This sphere has only one non-zero (reduced) Betti number, which is $\Betti{2\ell+1}{} = 1$.
  There are $\binom{k}{\ell+1}$ such full subcomplexes.
  The common intersection of any number of these subcomplexes is a complex of similar type, namely the full subcomplex of $\Delaunay{A}$ defined by the points on the common circles, which has the homotopy type of the $(2i+1)$-sphere, with $i \leq \ell$.
  By repeated application of the Mayer--Vietoris sequence \cite[page 149]{Hat02}, this implies that the Betti numbers of $\Alfa{\ell}{\ell}$ are bounded by a function of $k$ and are, thus, independent of $n$.
  Since we assume that $k$ is a constant, we have $\Betti{p}{(\Alfa{\ell}{\ell})} = O(1)$ for every $p$.
\end{proof}

Now we are ready to complete the proof of Theorem~\ref{thm:maximum_Betti_numbers_in_even_dimensions} for $p \leq 2k-2$.
To establish relation~\eqref{eqn:even1}, fix $p$ between $0$ and $k-1$ and consider $\Alfa{p}{-1} = \RadiusF^{-1} [0,\radius{}{p}{-1}]$, which is the Alpha complex consisting of all simplices that touch $p$ or fewer circles, together with all simplices that touch $p+1$ circles but each circle in only one point.
In other words, $\Alfa{p}{-1}$ is $\Alfa{p-1}{p-1}$ together with all the $\binom{k}{p+1} n^{p+1}$ $p$-simplices that have no short edges.
By Lemma~\ref{lem:constant_homology_in_even_dimensions}, $\Alfa{p-1}{p-1}$ has only a constant number of $(p-1)$-cycles.
Hence, only a constant number of the $p$-simplices can give death to $(p-1)$-cycles, while the remaining $p$-simplices give birth to $p$-cycles.
This is because every $p$-simplex either gives birth or death, so if it cannot give death to a $(p-1)$-cycle, then it gives birth to a $p$-cycle.
Hence, $\Betti{p}{(\Alfa{p}{-1})} = \binom{k}{p+1} n^{p+1} \pm O(1)$, as claimed.
The proof of relation~\eqref{eqn:even2} is similar but inductive. 
The induction hypothesis is
\begin{align}
    \Betti{p}{(\Alfa{k-1}{p-k})}  &=  \tbinom{k-1}{p-k+1} \cdot n^k \pm O(1) .
\end{align}
For $p = k-1$, it claims $\Betti{k-1}{(\Alfa{k-1}{-1})} = n^k \pm O(1)$, which is what we just proved.
In other words, relation~\eqref{eqn:even1} furnishes the base case at $p = k-1$.
A single inductive step takes us from $\Alfa{k-1}{p-k}$ to $\Alfa{k-1}{p-k+1}$; that is: we add all simplices that touch all $k$ circles and $p-k+2$ of them in two vertices to $\Alfa{k-1}{p-k}$.
The number of such simplices is the number of ways we can pick a pair of consecutive vertices from $p-k+2$ circles and a single vertex from the remaining $2k-p-2$ circles.
Since there are equally many vertices as there are consecutive pairs, this number is $\binom{k}{p-k+2} n^k$.
The dimension of these simplices is $(k-1)+(p-k+1)+1 = p+1$.
Some of these $(p+1)$-simplices give death to $p$-cycles, while the others give birth to $(p+1)$-cycles in $\Alfa{k-1}{p-k+1}$.
By the induction hypothesis, there are $\binom{k-1}{p-k+1} \cdot n^k \pm O(1)$ $p$-cycles in $\Alfa{k-1}{p-k}$, so this is also the number of $(p+1)$-simplices that give death.
Since $\binom{k}{p-k+2} - \binom{k-1}{p-k+1} = \binom{k-1}{p-k+2}$, this implies
\begin{align}
  \Betti{p}{(\Alfa{k-1}{p-k+1})}  &=  \tbinom{k-1}{p-k+2} \cdot n^k \pm O(1) ,
\end{align}
as required to finish the inductive argument.

\section{Three Dimensions}
\label{sec:3}

In this section, we answer the maximum Betti number question for \v{C}ech complexes in the smallest odd dimension in which it is non-trivial:
\begin{theorem}[Maximum Betti Numbers in $\Rspace^3$]
  \label{thm:maximum_Betti_numbers_in_three_dimensions}
  For every $n \geq 2$, there exist $N = 2n+2$ points in $\Rspace^3$ such that the \v{C}ech complex for a radius has first Betti number $\Betti{1}{} = (n+1)^2 - 1$ and for another radius has second Betti number $\Betti{2}{} = n^2$.
\end{theorem}
The proof consists of four steps: the construction of the set in Section~\ref{sec:3.1}, the analysis of the circumradii in Section~\ref{sec:3.2}, the argument that all simplices in the Delaunay mosaic are critical in Section~\ref{sec:3.3}, and the final counting of the tunnels and voids in Section~\ref{sec:3.4}.

\subsection{Construction}
\label{sec:3.1}

\begin{figure}[hbt]
    \centering \vspace{0.1in}
    \resizebox{!}{2.6in}{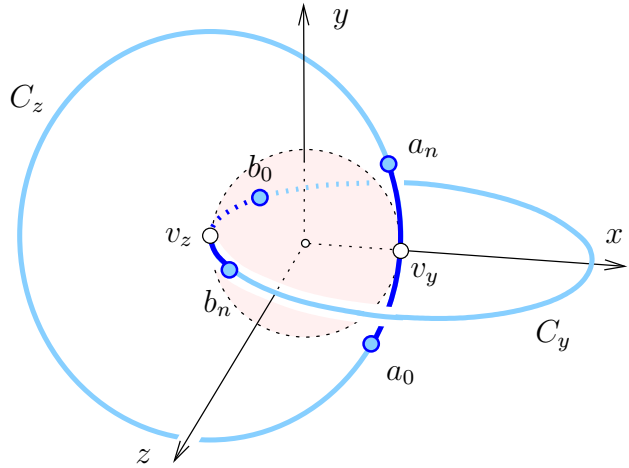}
    \vspace{-0.2in}
    \caption{Two linked unit circles in orthogonal coordinate planes of $\Rspace^3$, each touching the shaded sphere centered at the origin and each passing through the center of the other circle.
    There are $n+1$ points on each circle, on both sides and near the center of the other circle.}
    \label{fig:link}
\end{figure}

Given $n$ and $0 < \Delta < 1$, we construct the point set, $A = A_3 (n,\Delta)$, using two linked circles in $\Rspace^3$: $\Circle{z}$ with center $v_z = (-\frac{1}{2}, 0, 0)$ in the $xy$-plane defined by $(-\frac{1}{2} + \cos \varphi, \sin \varphi, 0)$ for $0 \leq \varphi < 2 \pi$, and $\Circle{y}$ with center $v_y = (\frac{1}{2}, 0, 0)$ in the $xz$-plane defined by $(\frac{1}{2} - \cos \psi, 0, \sin \psi)$ for $0 \leq \psi < 2 \pi$; see Figure~\ref{fig:link}.
On each circle, we choose $n+1$ points close to the center of the other circle.
To be specific, take the points $(0,-\Delta,0)$ and $(0,\Delta,0)$, and project them to $\Circle{z}$ along the $x$-axis. The resulting points are denoted by $a_0 = (-\frac{1}{2}+\sqrt{1-\Delta^2},-\Delta,0)$ and $a_n = (-\frac{1}{2}+\sqrt{1-\Delta^2},\Delta,0)$. Divide the arc between them into $n$ equal pieces by the points $a_1, a_2, \ldots, a_{n-1}$.
Symmetrically, project the points $(0,0,-\Delta)$ and $(0,0,\Delta)$ to $b_0 = (\frac{1}{2}-\sqrt{1-\Delta^2},0,-\Delta)$ and $b_n = (\frac{1}{2}-\sqrt{1-\Delta^2},0,\Delta)$ lying on $\Circle{y}$, and place $n-1$ points $b_1, b_2, \ldots, b_{n-1}$ on the arc between them, dividing it into $n$ equal pieces.
Let $\ee = \ee(n,\Delta)$ be the half-length of the (straight) edge connecting two consecutive points of either sequence.
Clearly, $\ee$ is a function of $n$ and $\Delta$, and it is easy to see that
\begin{align}
  \Delta / n < \ee < \tfrac{\pi}{2} \Delta / n
  \mbox{\rm ~~~and~~~} \ee \stackrel{\Delta \to 0}{\longrightarrow} \Delta/n .
  \label{eqn:epsilon}
\end{align}

A sphere that does not contain a circle intersects it in at most two points.
It follows that the sphere that passes through four points of $A$ is empty if and only if two of the four points are consecutive on one circle and the other two are consecutive on the other.
This determines the Delaunay mosaic: its $N = 2n+2$ vertices are the points $a_i$ and $b_j$, its $2n + (n+1)^2$ edges are of the forms $a_i a_{i+1}$, $b_j b_{j+1}$, and $a_i b_j$, its $2n(n+1)$ triangles are of the forms $a_i a_{i+1} b_j$ and $a_i b_j b_{j+1}$, and its $n^2$ tetrahedra of the form $a_i a_{i+1} b_j b_{j+1}$.
Keeping with the terminology introduced in Section~\ref{sec:2}, we call the edges $a_ib_j$ \emph{long} and the edges $a_ia_{i+1}$ and $b_j b_{j+1}$ \emph{short}.
Hence, every triangle in the Delaunay mosaic has one short and two long edges, and every tetrahedron has two short and four long edges.

\subsection{Divergence from the Ideal}
\label{sec:3.2}

The simplices in $\Delaunay{A}$ are not quite ideal, in the sense of Section~\ref{sec:2}.
We, therefore, need upper and lower bounds on their sizes, as quantified by their circumradii.
We will make repeated use of the following two inequalities, which both hold for $x > -1$:
\begin{align}
  \sqrt{1+x}  &\leq  1 + \tfrac{x}{2} , 
    \label{eqn:ineq1} \\
  \sqrt{1+x}  &\geq  1 + \tfrac{x}{2+x} .
    \label{eqn:ineq2}
\end{align}
For example, we will obtain some bounds on the radii of the triangle and tetrahedron in Lemma~\ref{lem:ideal_triangle_and_tetrahedron}, avoiding the use of square roots.
For the triangle, we rewrite \eqref{eqn:idealtriangle} to $4r^2(s) = 1 + x$ with $x = {s^2}/{(1-s^2)}$, and for the tetrahedron, we have $4R^2(s) = 1+x$ with $x = 2s^2$:
\begin{align}
  1 + \frac{1}{2} s^2  <  1 + \frac{s^2 / (1-s^2)}{2 + s^2 / (1-s^2)} \leq  \,2r(s)\,  &\leq  1 + \frac{s^2}{2-2s^2}  <  1 + \frac{10}{19} s^2 ,
    \label{eqn:triangleradius} \\
  1 + \frac{10}{11} s^2 \leq 1 + \frac{s^2}{1+s^2} \leq 2R(s)  &\leq 1+s^2 ,
    \label{eqn:tetrahedronradius}
\end{align}
where we assume that $n$ is large enough to imply $2 - 2 s^2 > 1.9$ and therefore $1+s^2 < 1.1$.
We begin by proving bounds on the lengths of long edges.
\begin{lemma}[Bounds for Long Edges in $\Rspace^3$]
  \label{lem:bounds_for_long_edges_in_three_dimensions}
  Let $0 < \Delta < 1$ and $A = A_3 (n, \Delta) \subseteq \Rspace^3$.
  Then the half-length of any long edge, $E \in \Delaunay{A}$, satisfies $\frac{1}{2} \leq R_E \leq \frac{1}{2} (1 + \Delta^4)$.
\end{lemma}
\begin{proof}
  To verify the lower bound, let $a \in \Circle{z}$ and consider the sphere with unit radius centered at $a$.
  This sphere intersects the $xz$-plane in a circle of radius at most $1$, whose center lies on the $x$-axis.
  The circle passes through $v_z \in \Circle{y}$, which implies that the rest of $\Circle{y}$ lies on or outside the circle and, therefore, on or outside the sphere centered at $a$.
  Hence, $\Edist{a}{b} \geq 1$ for all $b \in \Circle{y}$, which implies the required lower bound.

  \smallskip
  To establish the upper bound, observe that the distance between $a$ and $b$ is maximized if the two points are chosen as far as possible from the $x$-axis, so $4 R_E^2 \leq \Edist{a_0}{b_0}^2$.
  By construction, $a_0 = (-\frac{1}{2}+\sqrt{1-\Delta^2}, - \Delta, 0)$ and $b_0 = ( \frac{1}{2}-\sqrt{1-\Delta^2}, 0, -\Delta)$.
  Hence,
  \begin{align}
    4 R_E^2  &\leq  \norm{(-1+2\sqrt{1-\Delta^2}, -\Delta, \Delta)}^2
     =  5 - 2\Delta^2 - 4\sqrt{1-\Delta^2} 
       \label{eqn:long1} \\
    &\leq  5 - 2\Delta^2 - 4\left( 1 - \frac{\Delta^2}{2-\Delta^2} \right)
     =  1 + \frac{2\Delta^4}{2-\Delta^2}
       \label{eqn:long2} \\
    &\leq  1 + 2 \Delta^4 ,
     \label{eqn:long3}
  \end{align}
  where we used \eqref{eqn:ineq2} to get \eqref{eqn:long2} from \eqref{eqn:long1}, and $\Delta^2 < 1$ to obtain the final bound.
  Applying \eqref{eqn:ineq1}, wet get $2R_E \leq 1 + \Delta^4$, as required.
\end{proof}

Next, we estimate the circumradii of the triangles in $\Delaunay{A}$.
To avoid the computation of a constant, we use the big-Oh notation for $\Delta$, in which we assume that $n$ is fixed.
\begin{lemma}[Bounds for Triangles in $\Rspace^3$]
  \label{lem:bounds_for_triangles_in_three_dimensions}
  Let $0 < \Delta < \sqrt{2}/n$, $A = A_3 (n,\Delta) \subseteq \Rspace^3$, and $\ee = \ee(n,\Delta)$.
  Then the circumradius of any triangle, $F$, satisfies
  $\frac{1}{2} + \frac{1}{4} \ee^2  \leq  R_F  \leq  \frac{1}{2} + \frac{1}{4} \ee^2 +O(\Delta^4).$
\end{lemma}
\begin{proof}
  To see the lower bound, recall that the short edge of $F$ has length $2 \ee$ and the two long edges have lengths at least $1$.
  We place the endpoints of the short edge on a circle of radius $r(\ee)$.
  By the choice of the radius, there is only one point on this circle with distance at least $1$ from both endpoints, and it has distance $1$ from both.
  For any radius smaller than $r(\ee)$, there is no such point, which implies that the circumradius of $F$ satisfies $R_F \geq r(\ee) \geq \frac{1}{2} + \frac{1}{4} \ee^2$, where the second inequality follows from \eqref{eqn:triangleradius}.

  \smallskip
  To prove the upper bound, we draw $F$ in the plane, assuming its circumcircle is the circle with radius $R_F$ centered at the origin.
  Let $a, b, c$ be the vertices of $F$, where $a$ and $c$ are the endpoints of the short edge.
  We have $0 \in F$, since otherwise one of the angles at $a$ and $c$ is obtuse, in which case the squared lengths of the two long edges differ by at least $4 \ee^2$. 
  By assumption, $\sqrt{2} \Delta^2 < 2\Delta/n \leq 2 \ee$, in which we get the second inequality from \eqref{eqn:epsilon}.
  But this implies that the difference between the squared lengths of the two long edges is larger than $2 \Delta^4$, which contradicts \eqref{eqn:long3}.
  Hence, $b$ lies between the antipodes of the other two vertices, $a' = -a$ and $c' = -c$.
  By construction, $\Edist{a'}{c'} = 2\ee$.
  Assuming $\Edist{b}{a'} \leq \Edist{b}{c'}$, this implies
  \begin{align}
    \Edist{b}{a'}  &\leq  R_F \arcsin \tfrac{\ee}{R_F}
                    \leq  \arcsin \ee
                    =  \ee + O(\ee^3).
        \label{eqn:btoaprime}
  \end{align}
  Here, the second inequality follows from $R_F \geq 1$, using the convexity of the arcsin function, and the final expression using the Taylor expansion $\arcsin x = x + \frac{1}{6} x^3 + \frac{3}{40} x^5 + \ldots$.
  Now consider the triangle with vertices $a, a', b$.
  By the Pythagorean theorem,
  \begin{align}
    4 R_F^2  &=  \Edist{b}{a}^2 + \Edist{b}{a'}^2
            <  1 + 4 \Delta^4 + \ee^2 + O(\ee^4)
            =  1 + \ee^2 + O(\Delta^4) ,
  \end{align}
  where we used Lemma~\ref{lem:bounds_for_long_edges_in_three_dimensions} and \eqref{eqn:btoaprime} to bound $\Edist{b}{a}^2$ and $\Edist{b}{a'}^2$, respectively.
  We get the final expression using $\ee < \Delta$.
  Applying \eqref{eqn:ineq1}, we obtain $2 R_F \leq 1 + \frac{1}{2} \ee^2 + O(\Delta^4)$, as claimed.
\end{proof}

Similar to the case of triangles, it is not difficult to establish that the circumradius of any tetrahedron in the Delaunay mosaic is at least the circumradius of the ideal tetrahedron.
\begin{lemma}[Lower Bound for Tetrahedra in $\Rspace^3$]
  \label{lem:lower_bound_for_tetrahedra_in_three_dimensions}
  Let $0 < \Delta < 1$, $A = A_3 (n,\Delta) \subseteq \Rspace^3$, and $\ee = \ee(n,\Delta)$.
  Then the circumradius of any tetrahedron $T \in \Delaunay{A}$ satisfies $\frac{1}{2} + \frac{5}{11} \ee^2 \leq R_T$.
\end{lemma}
\begin{proof}
  By construction, $T$ has two disjoint short edges, both of length $2 \ee$.
  We place the endpoints of one short edge on a sphere of radius $R(\ee)$.
  The set of points on this sphere that are at distance at least $1$ from both endpoints is the intersection of two spherical caps whose centers are antipodal to the endpoints.
  We call this intersection a \emph{spherical bi-gon}.
  Since the two caps have the same size, the two corners of the bi-gon are further apart than any other two points of the bi-gon.
  By choice of the radius, $R(\ee)$, the edge connecting the two corners has length $2 \ee$.
  Hence, these corners are the only possible choice for the remaining two vertices of $T$, and for a radius smaller than $R(\ee)$, there is no choice.
  It follows that the circumradius of $T$ is at least $R(\ee)$, and we get the claimed lower bound from \eqref{eqn:tetrahedronradius}.
\end{proof}

\subsection{All Simplices are Critical}
\label{sec:3.3}

Since no empty sphere passes through more than four points of $A$, the Delaunay mosaic of $A$ is simplicial, and the radius function on this Delaunay mosaic is a generalized discrete Morse function~\cite{BaEd17}.
Furthermore, all simplices are critical; see Definition~\ref{def:critical_cell}.
The point set depends on two parameters, $n$ and $\Delta$, and we consider $n$ fixed while $\Delta$ goes to zero.
\begin{lemma}[All Critical in $\Rspace^3$]
  \label{lem:all_critical_in_three_dimensions}
  Let $n \geq 2$, $\Delta > 0$ sufficiently small, and $A = A_3 (n,\Delta) \subseteq \Rspace^3$.
   Then every simplex of the Delaunay mosaic of $A$ is critical.
\end{lemma}
\begin{proof}
  It is clear that the vertices and the short edges are critical, but the other simplices in $\Delaunay{A}$ require an argument.
  We begin with the long edges.
  Fix $i$ and $j$, and write $\Sphere{2} (i; j)$ for the smallest sphere that passes through $a_i$ and $b_j$.
  Its center is the midpoint of the long edge and by \eqref{eqn:long3} its squared diameter is between $1$ and $1 + 2 \Delta^4$.
  The distance between $a_i$ and any $a_\ell$, $\ell \neq i$, is at least $2 \ee$.
  Assuming $a_\ell$ is on or inside $\Sphere{2} (i;j)$, we thus have $\Edist{a_\ell}{b_j}^2 \leq 1 + 2 \Delta^4 - 4 \ee^2$, which, for sufficiently small $\Delta > 0$, is less than $1$.
  But this contradicts the lower bound in Lemma~\ref{lem:bounds_for_long_edges_in_three_dimensions}, so $a_\ell$ lies outside $\Sphere{2} (i;j)$.
  By a symmetric argument, all $b_\ell$, $\ell \neq j$, lie outside $\Sphere{2} (i;j)$.
  Hence, $\Sphere{2} (i;j)$ is strictly empty, for all $0 \leq i, j \leq n$, which implies that all edges of $\Delaunay{A}$ are critical edges of the radius function.
  
  \smallskip
  The fact that all edges of $\Delaunay{A}$ are critical implies that all triangles are acute.
  Indeed, if $a_i b_j b_{j+1}$ is not acute, then the midpoint of one long edge is at least as close to the third vertex as to the endpoints of the edge.
  Hence, any non-acute triangle would be an obstacle to the criticality of an edge, which implies that no such triangle can exist.
  However, the fact that all triangles are acute does not imply that all of them are critical.
  To prove the criticality of the Delaunay triangles, let $x$ be the circumcenter of $a_i b_j b_{j+1}$, let $\Sphere{2} (i; j, j+1)$ be centered at $x$ and pass through $a_i$, $b_j$, $b_{j+1}$, and let $a$ be the point other than $a_i$ in which $\Sphere{2} (i;j,j+1)$ intersects $\Circle{z}$.
  Since $a_i b_j b_{j+1}$ is acute, $x$ lies in the interior of the triangle. 
  It remains to show that the sphere is strictly empty.
  To this end, let $x'$ and $x''$ be the centers of $\Sphere{2} (i;j)$ and $\Sphere{2} (i;j+1)$, let $a'$ and $a''$ be the points other than $a_i$ in which the two spheres intersect $\Circle{z}$, and consider the lines that pass through $x$ and $x'$ and through $x$ and $x''$, respectively.
  Note that $x$ lies between $x'$ and $x''$.
  This implies that $a$ is between $a'$ and $a''$.
  Since $\Sphere{2} (i;j)$ and $\Sphere{2} (i;j+1)$ are strictly empty, $a'$ and $a''$ lie strictly between $a_{i-1}$ and $a_{i+1}$, and so does $a$.
  Hence, $\Sphere{2} (i; j, j+1)$ is strictly empty, which implies that all triangles of $\Delaunay{A}$ are critical triangles of the radius function.

  \smallskip
  Since all triangles are critical, all tetrahedra of $\Delaunay{A}$ must also be critical. One can argue in two ways. Combinatorially: the radius function pairs non-critical tetrahedra with non-critical triangles, but there are no such triangles. Geometrically: since every triangle has a non-empty intersection with its dual Voronoi edge, every tetrahedron must contain its dual Voronoi vertex. 
\end{proof}

\subsection{Counting the Tunnels and Voids}
\label{sec:3.4}

Before counting the tunnels and voids, we recall that $\RadiusF \colon \Delaunay{A} \to \Rspace$ maps each simplex to the radius of its smallest empty sphere that passes through its vertices.
By Lemma~\ref{lem:all_critical_in_three_dimensions}, all simplices of $\Delaunay{A}$ are critical, so $\RadiusF (E)$ is equal to the circumradius of $E$, for every edge $E \in \Delaunay{A}$, and similarly for every triangle and every tetrahedron.
\begin{corollary}[Ordering of Radii in $\Rspace^3$]
  \label{cor:ordering_of_radii_in_three_dimensions}
  Let $\Delta > 0$ be sufficiently small, let $A = A_3 (n,\Delta) \subseteq \Rspace^3$, and let $\RadiusF \colon \Delaunay{A} \to \Rspace$ be the radius function.
  Then $\RadiusF (E) < \RadiusF (F) < \RadiusF (T)$ for every edge $E$, triangle $F$, and tetrahedron $T$ in $\Delaunay{A}$.
\end{corollary}
\begin{proof}
  Using Lemma~\ref{lem:bounds_for_long_edges_in_three_dimensions} for the edges,  Lemma~\ref{lem:bounds_for_triangles_in_three_dimensions} for the triangles, and 
  Lemma~\ref{lem:lower_bound_for_tetrahedra_in_three_dimensions} for the tetrahedra in the Delaunay mosaic of $A$, we get
  \begin{align}
    \RadiusF (E) = R_E  &<     \tfrac{1}{2} + O(\Delta^4), 
        \label{eqn:edge} \\
    \tfrac{1}{2} + \tfrac{1}{4} \ee^2  \leq  \RadiusF (F) = R_F  &<     \tfrac{1}{2} + \tfrac{1}{4} \ee^2 + O(\Delta^4),
        \label{eqn:triangle} \\
    \tfrac{1}{2} + \tfrac{5}{11} \ee^2   \leq  \RadiusF (T) = R_T  & ,
        \label{eqn:tetrahedron}
  \end{align}
  so for sufficiently small $\Delta > 0$, the edges precede the triangles, and the triangles precede the tetrahedra in the filtration of the simplices.
\end{proof}

For the final counting, choose $\rho_1$ to be any number strictly between the maximum radius of any edge and the minimum radius of any triangle. The existence of such a number is guaranteed by Corollary~\ref{cor:ordering_of_radii_in_three_dimensions}.
The corresponding \v{C}ech complex is the $1$-skeleton of the Delaunay mosaic.
It is connected, with $N = 2n+2$ vertices and $2n + (n+1)^2$ edges.
The number of independent cycles is the difference plus $1$, which implies $\Betti{1}{(\Cech{A}{\rho_1})} = (n+1)^2 - 1$, as claimed.
Similarly, choose $\rho_2$ between the maximum radius of any triangle and the minimum radius of any tetrahedron, which is again possible, by Corollary~\ref{cor:ordering_of_radii_in_three_dimensions}.
The corresponding \v{C}ech complex is the $2$-skeleton of the Delaunay mosaic.
The number of independent $2$-cycles is the number of missing tetrahedra.
This implies $\Betti{2}{(\Cech{A}{\rho_2})} = n^2$, as claimed.

\section{Odd Dimensions}
\label{sec:4}

In this section, we generalize the $3$-dimensional results presented in Section~\ref{sec:3} to every odd dimension.
\begin{theorem}[Maximum Betti Numbers in $\Rspace^{2k+1}$]
  \label{thm:maximum_Betti_numbers_in_odd_dimensions}
  For every $d = 2k+1 \geq 1$, $n \geq 2$, and sufficiently small $\Delta > 0$, there are a set $A = A_d(n,\Delta) \subseteq \Rspace^{2k+1}$ of $N = (k+1)(n+1)$ points and radii $\rho_0 < \rho_1 < \ldots < \rho_{2k}$ such that
  \begin{align}
    \Betti{p}{(\Cech{A}{\rho_p})} &= \tbinom{k+1}{p+1} \cdot (n+1)^{p+1} \pm O(1), \mbox{\rm ~~~~~for~} 0 \leq p \leq k ; 
      \label{eqn:thm41A} \\
    \Betti{p}{(\Cech{A}{\rho_p})} &= \tbinom{k}{p-k} \cdot (n+1)^{k+1} \pm O(n^k), \mbox{\rm ~~~for~} k+1 \leq p \leq 2k .
      \label{eqn:thm41B}
  \end{align}
\end{theorem}
The steps in the proof are the same as in Sections~\ref{sec:2} and \ref{sec:3}:  construction of the points, analysis of the circumradii, argument that all simplices are critical, and final counting of the cycles.
In contrast to the earlier sections, the analytic part of the proof is inductive and distinguishes between erecting a pyramid or a bi-pyramid on top of a lower-dimensional simplex.

\subsection{Construction}
\label{sec:4.1}

\begin{figure}[hbt]
  \centering \vspace{0.0in}
  \resizebox{!}{2.5in}{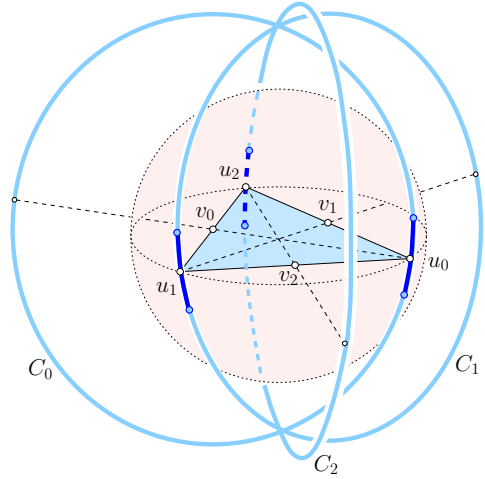}
  \caption{The projection of the $5$-dimensional construction to $\Rspace^3$, in which $x_3, x_4, x_5$ are all mapped to the same, vertical coordinate direction.
  The circles $\Circle{0}, \Circle{1}, \Circle{2}$ touch the shaded sphere in the vertices of the triangle.
  In $\Rspace^5$, the three circles belong to mutually orthogonal $2$-planes, so the two common points of the three circles in the drawing are an artifact of the particular projection.}
  \label{fig:cartoon}
\end{figure}

Equip $\Rspace^d$ with Cartesian coordinates, $x_1, x_2, \ldots, x_d$, and consider a regular $k$-simplex, denoted by $\Simplex{}$, in the $k$-plane spanned by $x_1, x_2, \ldots, x_k$.
It is not important where $\Simplex{}$ is located inside the coordinate $k$-plane, but we assume for convenience that its barycenter is the origin of the coordinate system.
It is, however, important that all edges of $\Sigma$ have unit length.
We will repeatedly need the squared circumradius, height, and in-radius of $\Sigma$, for which we state simple formulas for later convenience:
\begin{align}
  R_k^2 = \tfrac{k}{2(k+1)}, ~~~
  H_k^2 = \tfrac{k+1}{2k}, ~~~
  D_k^2 = R_k^2 - R_{k-1}^2 = \tfrac{1}{2 k (k+1)} .
    \label{eqn:regular_simplex}
\end{align}
Observe that the angle, $\alpha$, between an edge and a height of $\Simplex{}$ that meet at a shared vertex satisfies $\cos \alpha = H_k$.
Let $u_0, u_1, \ldots, u_k$ be the vertices of $\Simplex{}$, and let $v_\ell$ be the barycenter of the $(k-1)$-face opposite to $u_\ell$.
For each $0 \leq \ell \leq k$, consider the $2$-plane spanned by $u_\ell - v_\ell$ and the $x_{k+\ell+1}$-axis, and let $\Circle{\ell}$ be the circle in this $2$-plane, centered at $v_\ell$, that passes through $u_\ell$; see Figure~\ref{fig:cartoon}.
Its radius is the height of the $k$-simplex: $\radiusC = H_k$.
Given a global choice of the parameter, $0 < \Delta < H_k$, we cut $\Circle{\ell}$ at $x_{k+\ell+1} = \pm \Delta$ into four arcs and place $n+1$ point at equal angles along the arc that passes through $u_\ell$.
Repeating this step for each $\ell$, we get a set of $N = (k+1)(n+1)$ points, denoted $A = A_{2k+1} (n,\Delta)$.

A $(d-1)$-sphere that contains none of the circles $\Circle{\ell}$ intersects the $k+1$ circles in at most two points each.
It follows that a sphere that passes through $2k+2$ points of $A_d$ is empty if and only if it passes through two consecutive points on each of the $k+1$ circles.
This determines the Delaunay mosaic, which consists of $n^{k+1}$ $d$-simplices together with all their faces.
It follows that the number of $p$-simplices in $\Delaunay{A}$ is at most some constant times $n^m$, in which $m = \min \{ p+1, k+1 \}$ and the constant depends on $d = 2k+1$.
Building on the notation introduced in Section~\ref{sec:2}, we describe a simplex, $S \in \Delaunay{A}$, with two integers:  $\ell = \touch{S}$ is one less than the number of circles it touches, and $j = \short{S}$ is one less than the number of short edges.
Hence, $p = \ell+j+1$ is the dimension.
For each $0 \leq p \leq k$, there are $\binom{k+1}{p+1} (n+1)^{p+1}$ $p$-simplices that touch $\ell+1 = p+1$ circles and thus have $j+1 = 0$ short edges.
As suggested by a comparison with relation \eqref{eqn:thm41A} in Theorem~\ref{thm:maximum_Betti_numbers_in_odd_dimensions}, these $p$-simplices will be found responsible for the $p$-cycles counted by the $p$-th Betti number.

\subsection{Inductive Analysis}
\label{sec:4.2}

The bulk of the proof of Theorem~\ref{thm:maximum_Betti_numbers_in_odd_dimensions} is devoted to the analysis of the Delaunay simplices.
The goal is to prove bounds on the circumradii that are strong enough to separate simplices of different types, and to show that all simplices are critical.
The analysis is inductive with three hypotheses: the first about the circumradius, the second about the circumcenter, and the third about the projection of a vertex onto the affine hulls of the opposite facet.
To formulate the second hypothesis, we write $D_S$ for the radius of the largest ball centered at the circumcenter that is contained in a simplex, $S$.
To formulate the third hypothesis, we call a point $x \in \aff{S}$ \emph{edge-centric} if the distance between the projection of $x$ onto any edge of $S$ has distance at most $X_E = n\Delta^3$ from the midpoint of that edge, and we write $X_S$ for the maximum distance between any edge-centric point and the circumcenter of $S$.
Recall that $\ee = \ee(n,\Delta)$ is a function of $n$ and $\Delta$ that satisfies $\Delta/n \leq \ee \leq \frac{\pi}{2} \Delta/n$.

\medskip \noindent \textbf{Hypothesis I:}
    $R_S^2 = R_\ell^2 + \frac{j+1}{(\ell+1)^2} \ee^2 \pm O(\ee^3)$.
    
\smallskip \noindent \textbf{Hypothesis II:}
    $D_S^2 = \left\{ \begin{array}{cl}
                     D_\ell^2 \pm O(\ee^2) & \mbox{\rm if~} \short{S} = -1, \\
                     \frac{1}{(\ell+1)^2} \ee^2 \pm O(\ee^3) & \mbox{\rm if~} 0 \leq \short{S} \leq \touch{S};
                   \end{array} \right.$

\smallskip \noindent \textbf{Hypothesis III:}
    $X_S = O(\Delta^3)$,

\medskip \noindent
in which the big-Oh notation is used to suppress multiplicative constants, as usual.
We assume that $\Delta$ is chosen independent of the number of points, so in this context, $n$ is considered to be a constant, and we write $\Delta = O(\ee)$, for example.
The base case for the first two hypotheses will be covered by Lemmas~\ref{lem:bounds_for_long_edges_in_odd_dimensions} and \ref{lem:bounds_for_bisectors_in_odd_dimensions}, and
the third hypothesis holds for edges, by definition.
We will distinguish between two kinds of inductive steps, one reasoning from $(\ell-1,j)$ to $(\ell,j)$ and the other from $(\ell,j-1)$ to $(\ell,j)$.
We need some notions to describe the difference.
A \emph{facet} of a simplex is a face whose dimension is $1$ less than that of the simplex.
We call a vertex $a$ of $S$ a \emph{twin} if it is the endpoint of a short edge, in which case we write $a''$ for the other endpoint of that edge.
If $a$ is not a twin, we write $Q = S-a$ for the opposite facet, and call the pair $(a, Q)$ a \emph{pyramid} with \emph{apex} $a$ and \emph{base} $Q$.
The point of Hypothesis~III is that together with Lemma~\ref{lem:bounds_for_bisectors_in_odd_dimensions}, it will imply that $a$ projects to a point in $Q$ whose distance from the circumcenter of $Q$ is at most $X_S$.
If $a$ is a twin, then there are two pyramids, $(a,P)$ and $(a'',P)$ with $P=S-a-a''$, and we call this the \emph{bi-pyramid case}.

\subsubsection{Base Case}
\label{sec:4.2.1}

The only non-trivial base cases are when $S$ is a long edge in Hypothesis~I, and when $S$ is a short edge in Hypothesis~III.
To prove bounds on the length of a long edge, we write $R_E$ for its half-length, which is also its circumradius.
\begin{lemma}[Bounds for Long Edges in $\Rspace^{2k+1}$]
  \label{lem:bounds_for_long_edges_in_odd_dimensions}
  Let $d = 2k+1$, $0 < \Delta < 1$, and $A = A_d (n,\Delta) \subseteq \Rspace^d$.
  Then the squared length of any long edge satisfies $1 \leq 4 R_E^2 \leq 1 + 2 \Delta^4$.
\end{lemma}
\begin{proof}
  We simplify the computations by assuming that the endpoints $a$ and $b$ of $E$ are at equal distance from $\aff{\Simplex{}}$.
  Call this distance $\Delta$, suppose $a \in \Circle{0}$ and $b \in \Circle{1}$, and write $a'$ and $b'$ for their projections onto $\aff{\Simplex{}}$.
  Recall that $u_0$ is the point shared by $\Simplex{}$ and $\Circle{0}$, and note that $\Edist{a'}{u_0} = \xi = \radiusC - \sqrt{\radiusC^2 - \Delta^2}$, in which $\radiusC$ is the radius of $C_0$.
  Similarly, $\Edist{b'}{u_1} = \xi$.
  Let $\alpha$ be the angle enclosed by an edge of $\Simplex{}$ and a height of $\Simplex{}$ that shares a vertex with the edge.
  Set $\eta = \xi \cos \alpha$ and note that $\Edist{a'}{b'} = 1-2\eta$.
  By construction of $\Simplex{}$ as a regular simplex with unit length edges, we have $\cos \alpha = \radiusC$, so
  \begin{align}
    \Edist{a}{b}^2  &=  (1-2\eta)^2 + \Delta^2 + \Delta^2
       =  \left( 1 - 2 \radiusC^2 + 2 \radiusC \sqrt{\radiusC^2 - \Delta^2} \right)^2 + 2 \Delta^2 \\
      &=  \left( 1-2\radiusC^2 \right)^2 + 4\radiusC^2 \left( \radiusC^2 - \Delta^2 \right) + \left( 2-4\radiusC^2 \right) 2 \radiusC \sqrt{\radiusC^2 - \Delta^2} + 2 \Delta^2 \\
      &=  \left( 1-4\radiusC^2+8\radiusC^4 \right) - \left( 4\radiusC^2 - 2 \right) \left[ \Delta^2 + 2\radiusC \sqrt{\radiusC^2-\Delta^2} \right ].
        \label{eqn:longedgeodd}
  \end{align}
  The squared radius of the circles is $\radiusC^2 = (k+1)/(2k) > \frac{1}{2}$, which implies $4\radiusC^2 - 2 > 0$.
  Hence, we can bound $\Edist{a}{b}^2$ from below using \eqref{eqn:ineq1} to get $\sqrt{\radiusC^2-\Delta^2} \leq \radiusC \left[ 1 - \Delta^2/(2 \radiusC^2) \right]$.
  Plugging this inequality into \eqref{eqn:longedgeodd} and applying a sequence of elementary algebraic manipulations gives $\Edist{a}{b}^2 \geq 1$, as claimed.
  To prove the upper bound, we use \eqref{eqn:ineq2} to get $\sqrt{\radiusC^2-\Delta^2} \geq \radiusC \left[ 1 - \Delta^2/(2\radiusC^2-\Delta^2) \right]$.
  Plugging this inequality into \eqref{eqn:longedgeodd} gives
  \begin{align}
    \Edist{a}{b}^2  &\leq  \left( 1-4\radiusC^2+8\radiusC^4 \right) - \left( 4\radiusC^2-2 \right) \left[ \Delta^2+2\radiusC^2- \frac{2 \radiusC^2 \Delta^2}{2 \radiusC^2 - \Delta^2} \right] \\
      &=  1 + \left( 4 \radiusC^2 - 2 \right) \frac{\Delta^4}{2 \radiusC^2 - \Delta^2}
       \leq  1 + 2 \Delta^4 ,
  \end{align}
  where, to get the final inequality, we used that $\Delta^2 < 1$.
\end{proof}

If we first take the square root and then divide by $2$, we get $R_E \leq \frac{1}{2} (1 + \Delta^4)$ for the half-length or circumradius of the edge.
Since the length of long edges is so tightly controlled, the triangles formed by three long edges are almost equilateral, and the triangles formed by one short and two long edges are almost isosceles.
The next lemma quantifies this claim.
\begin{lemma}[Bounds for Bisectors in $\Rspace^{2k+1}$]
  \label{lem:bounds_for_bisectors_in_odd_dimensions}
  Let $d = 2k+1$, $\Delta > 0$ sufficiently small, and $A = A_d (n,\Delta) \subseteq \Rspace^d$.
  Then the distance between a vertex connected by long edges to the endpoints of another (short or long) edge and the bisector of this edge is at most $n \Delta^3 / 2$.
\end{lemma}
\begin{proof}
  Consider a vertex, $a$, connected by long edges to the endpoints, $b$ and $c$, of another (short or long) edge.
  Let $\delta$ be the distance of $a$ from the bisector of $b$ and $c$, which is maximized if the length difference is as large as possible while $\Edist{b}{c}$ is as small as possible.
  In this case, Pythagoras' theorem implies $(1+2\Delta^4) - (\ee+\delta)^2 = 1 - (\ee-\delta)^2$.
  Canceling $1$, $\ee^2$, and $\delta^2$ on both sides, we get $\Delta^4 = 2 \ee \delta$.
  Since $n \ee \geq \Delta$, this implies that $\delta = \Delta^4 / (2 \ee) \leq n \Delta^3 / 2$.
\end{proof}

We mention that choosing $\Delta$ is independent of $n$, so in this context, $n$ is considered a constant and we write $n \Delta^3 = O(\Delta^3)$.
We also note that the upper bound on the distance of a point connected by two long edges to the endpoints of a short edge from the bisector of these two points can be improved to $2 \Delta$.
We prefer the weaker bound in Lemma~\ref{lem:bounds_for_bisectors_in_odd_dimensions} because of its elementary proof.

\subsubsection{Inductive Step (Pyramid Case)}
\label{sec:4.2.2}

The inductive step consists of two lemmas.
The first one justifies the first kind of inductive step, from $(\ell-1,j)$ to $(\ell,j)$. It handles the transition from the base of a pyramid to the pyramid.
Letting $(a,Q)$ be a pyramid of $S$, we write $H_{Q,S}$ and $D_{Q,S}$ for the distances of $a$ and $\ccz{S}$ from $\aff{Q}$, respectively.
\begin{lemma}[Pyramid Step]
  \label{lem:pyramid_step}
  Let $d = 2k+1$, $\Delta > 0$ sufficiently small, $A = A_d (n,\Delta) \subseteq \Rspace^d$, and $\ee = \ee(n,\Delta)$.
  Furthermore, let $S \in \Delaunay{A}$, write $\ell = \touch{S}$ and $j = \short{S}$, assume $j < \ell$, and let $(a, Q)$ be a pyramid of $S$.
  Assuming $Q$ satisfies Hypotheses~I, II, and III, we have
  \begin{align}
    H_{Q,S}^2  &=  H_\ell^2 - \frac{j+1}{\ell^2} \ee^2 \pm O(\ee^3) ; 
      \label{eqn:L44A} \\
    D_{Q,S}^2  &=  D_\ell^2 - \frac{(2\ell+1)(j+1)}{\ell^2(\ell+1)^2} \ee^2 \pm O(\ee^3) ; 
      \label{eqn:L44B} \\
    R_S^2      &=  R_\ell^2 + \frac{j+1}{(\ell+1)^2} \ee^2 \pm O(\ee^3) ; 
      \label{eqn:L44C} \\
    X_S    &=  O(\Delta^3).
      \label{eqn:L44D}
  \end{align}
\end{lemma}
\begin{proof}
  By construction, $\touch{Q} = \ell-1$ and $\short{Q} = j$.
  Assume first that the projection of $a$ onto $\aff{Q}$ is $\ccz{Q}$.
  In this case, all edges connecting $a$ to $Q$ have the same length, $2 R_E$.
  Pythagoras' theorem implies $H_{Q,S}^2 = 4 R_E^2 - R_Q^2$.
  Using Lemma~\ref{lem:bounds_for_long_edges_in_odd_dimensions} and Hypothesis~I, we get the bounds for the squared height claimed in \eqref{eqn:L44A}:
  \begin{align}
    4 R_E^2  &=  1 \pm O(\Delta^4) ; 
      \label{eqn:Step1E} \\
    R_Q^2         &=  R_{\ell-1}^2 + \frac{j+1}{\ell^2} \ee^2 \pm O(\ee^3) ; 
      \label{eqn:Step1R} \\
    H_{Q,S}^2     &=  H_\ell^2 - \frac{j+1}{\ell^2} \ee^2 \pm O(\ee^3) ,
      \label{eqn:Step1H}
  \end{align}
  where \eqref{eqn:Step1H} follows from \eqref{eqn:Step1E} and \eqref{eqn:Step1R}, using $1 - R_{\ell-1}^2 = H_\ell^2$.
  This proves \eqref{eqn:L44A}.
  Since $(H_{Q,S} - D_{Q,S})^2 = R_S^2$ and $R_Q^2 + D_{Q,S}^2 = R_S^2$, we get $H_{Q,S}^2 - 2 D_{Q,S} H_{Q,S} = R_Q^2$.
  Therefore,
  \begin{align}
    D_{Q,S}  &=  \frac{H_{Q,S}^2 - R_Q^2}{2 H_{Q,S}}
              ~~=  \frac{1}{2} H_{Q,S} - \frac{1}{2} \frac{R_Q^2}{H_{Q,S}} ; 
      \label{eqn:D} \\
    R_S      &=  H_{Q,S} - D_{Q,S}
              =  \frac{1}{2} H_{Q,S} + \frac{1}{2} \frac{R_Q^2}{H_{Q,S}} .
      \label{eqn:R}
  \end{align}
  Using the formulas for $R_\ell$, $H_\ell$, $D_\ell$ in \eqref{eqn:regular_simplex}, it is easy to prove the corresponding relations for the regular $\ell$-simplex:
  $D_\ell = \frac{1}{2} H_\ell - \frac{1}{2} R_{\ell-1}^2 / H_\ell$ and $R_\ell = \frac{1}{2} H_\ell + \frac{1}{2} R_{\ell-1}^2 / H_\ell$.
  Starting with \eqref{eqn:Step1H}, we use $\sqrt{1-x} = 1 - \frac{x}{2} + \ldots$ and $1 / \sqrt{1-x} = 1 + \frac{x}{2} + \ldots$ to get
  \begin{align}
    H_{Q,S}  &=  H_\ell - \frac{j+1}{2 \ell^2 H_\ell} \ee^2 \pm O(\ee^3) ; 
      \label{eqn:H} \\
    \frac{1}{H_{Q,S}}  &=  \frac{1}{H_\ell} + \frac{j+1}{2 \ell^2 H_\ell^3} \ee^2 \pm O(\ee^3) ; 
      \label{eqn:1overH} \\
    \frac{R_Q^2}{H_{Q,S}}  &=  \frac{R_{\ell-1}^2}{H_\ell} + \left[ \frac{j+1}{\ell^2 H_\ell} + \frac{R_{\ell-1}^2 (j+1)}{2 \ell^2 H_\ell^3} \right] \ee^2 \pm O(\ee^3) .
      \label{eqn:RoverH}
  \end{align}
  We plug \eqref{eqn:H} and \eqref{eqn:RoverH} into \eqref{eqn:D} and \eqref{eqn:R}, while using the relations for $D_\ell$ and $R_\ell$ mentioned above, as well as ${R_\ell}/{H_\ell} = {\ell}/{(\ell+1)}$, ${D_\ell}/{H_\ell} = {1}/{(\ell+1)}$, and ${R_{\ell-1}^2}/{H_\ell^2} = {(\ell-1)}/{(\ell+1)}$:
  \begin{align}
    D_{Q,S}  &=  \left[ \frac{1}{2} H_\ell - \frac{1}{2} \frac{R_{\ell-1}^2}{H_\ell} \right] - \left[ \frac{j+1}{4 \ell^2 H_\ell} + \frac{j+1}{2 \ell^2 H_\ell} + \frac{R_{\ell-1}^2 (j+1)}{4 \ell^2 H_\ell^3} \right] \ee^2 \pm O(\ee^3) \nonumber \\
             &=  D_\ell - \frac{(2\ell+1)(j+1)}{2 \ell^2 (\ell+1)^2 D_\ell} \ee^2 \pm O(\ee^3) ; \\
    R_S      &=  \left[ \frac{1}{2} H_\ell + \frac{1}{2} \frac{R_{\ell-1}^2}{H_\ell} \right] + \left[ - \frac{j+1}{4 \ell^2 H_\ell} + \frac{j+1}{2 \ell^2 H_\ell} + \frac{R_{\ell-1}^2 (j+1)}{4 \ell^2 H_\ell^3} \right] \ee^2 \pm O(\ee^3) \nonumber \\
             &=  R_\ell + \frac{j+1}{2 (\ell+1)^2 R_\ell} \ee^2 \pm O(\ee^3) .
  \end{align}
  Taking squares, we get \eqref{eqn:L44B} and \eqref{eqn:L44C}, but mind that this is only for the special case in which the apex projects orthogonally to the circumcenter of the base.
  To prove the bounds in the general case, we recall that Hypothesis~III asserts that the projection of $a$ onto $\aff{Q}$ is at most $O(\Delta^3)$ units of length from $\ccz{Q}$.
  Hence, we get an additional error term of $O(\Delta^3)$ in all the above equations, but this does not change any of the bounds as stated.
  
  \smallskip
  It remains to prove \eqref{eqn:L44D}.
  By the inductive assumption, we have $X_Q = O(\Delta^3)$.
  Consider the locus of points in $\aff{S}$ whose projections to $\aff{Q}$ are at distance at most $X_Q$ from $\ccz{Q}$. This is a solid cylinder. 
  In addition, consider the locus of points whose projections to an edge connecting $a$ to a vertex of $Q$ are at distance at most $X_E$ from the midpoint of this edge. This is a slab between two parallel hyperplanes in $\aff{S}$.
  The points at distance at most $X_S$ from $\ccz{S}$ are contained in the intersection of this cylinder and the slab.
  Since $H_\ell^2 = {(\ell+1)}/{(2 \ell)}$ is strictly larger than $R_{\ell-1}^2 = {(\ell-1)}/{(2 \ell)}$, the angle at which the central axis of the cylinder and the central hyperplane of the slab intersect is larger than $\pi/4$, provided that $\Delta >0$ is sufficiently small.
  But then the intersection is contained in a ball of radius at most $\sqrt{2} X_Q + X_E = O(\Delta^3)$.
\end{proof}

Note that $D_S$ is the minimum of the $D_{Q,S}$, over all facets $Q$ of $S$.
Hence, Lemma~\ref{lem:pyramid_step} proves Hypothesis~II in the case in which $S$ has no short edges.

\subsubsection{Inductive Step (Bi-pyramid Case)}
\label{sec:4.2.3}

The second kind of inductive step---from $(\ell, j-1)$ to $(\ell, j)$---makes use of a distance function between affine subspaces of $\Rspace^d$.
Such a function is nonnegative, by definition, as well as convex; see e.g.\ Rockafellar~\cite[pages 28 and 34]{Roc70}.
In our case, the function will measure the distance from a $p$-plane to a $(d-1)$-plane, so it has a well-defined gradient, provided that the distance is taken with a sign, which is different on the two sides of the intersection with the hyperplane.
\begin{lemma}[Bi-pyramid Step]
  \label{lem:bi-pyramid_step}
  Let $d = 2k+1$, $\Delta > 0$ sufficiently small, $A = A_d (n,\Delta) \subseteq \Rspace^d$, and $\ee = \ee(n,\Delta)$.
  Furthermore, let $S \in \Delaunay{A}$, with $\ell = \touch{S}$ and $j = \short{S} \geq 0$, and let $a$ and $a''$ be the endpoints of a short edge.
  Assuming $Q = S-a''$ and $Q'' = S-a$ satisfy Hypotheses~I, II, and III, we have
  \begin{align}
    D_{Q,S}^2  &=  \frac{1}{(\ell+1)^2} \ee^2 \pm O(\ee^3); 
      \label{eqn:L45A} \\
    R_S^2      &=  R_\ell^2 + \frac{j+1}{(\ell+1)^2} \ee^2 \pm O(\ee^3); 
      \label{eqn:L45B} \\
    X_S        &=  O(\Delta^3).
      \label{eqn:L45C}
  \end{align}
\end{lemma}
\begin{proof}
  By construction, $\touch{Q} = \touch{Q''} = \ell$, $\short{Q} = \short{Q''} = j-1$, and $(a, Q-a)$ and $(a'', Q''-a'')$ are pyramids.
  We write $P = Q-a = Q''-a''$ for the common base, which has $\touch{P} = \ell-1$ and $\short{S} = j-1$..
  Let $M$ be the bisector of $a$ and $a''$.
  It intersects the short edge orthogonally at its midpoint.
  Writing $\psi \colon \aff{Q} \to M$ for the distance function from $\aff{Q}$ to $M$, we have $\psi(a) = \ee$ and, by Lemma~\ref{lem:bounds_for_bisectors_in_odd_dimensions}, $\psi(b) \leq n\Delta^3$, for all vertices $b$ of $P$.
  Let $a'$ be the projection of $a$ onto $\aff{P}$.
  By Hypotheses II and III, $a'$ is closer to $\ccz{P}$ than the radius of the largest ball centered at $\ccz{P}$ which is contained in $P$.
  Hence, $a' \in P$, so $\psi(a') \leq n\Delta^3$ by the convexity of the distance function.
  The signed version of $\psi$ is linear and, thus, has a well-defined gradient.
  To compute it, recall Lemma~\ref{lem:pyramid_step}, which shows that the height of $(a,P)$ and $\Edist{\ccz{Q}}{\ccz{P}}$ satisfy
  \begin{align}
    H_{P,Q}^2  &=  H_\ell^2 - \frac{j}{\ell^2} \ee^2 \pm O(\ee^3) ; 
      \label{eqn:L45height} \\
    D_{P,Q}^2  &=  D_\ell^2 - \frac{(2\ell+1) j}{\ell^2 (\ell+1)^2} \ee^2 \pm O(\ee^3) .
      \label{eqn:L45depth}
  \end{align}
  By \eqref{eqn:L45height}, the gradient of $\psi$ has length $\norm{\nabla \psi} = \ee / H_{P,Q} \pm O(\Delta^3) = \ee / H_\ell \pm O(\ee^3)$, and by \eqref{eqn:L45depth}, the value of the function at the circumcenter is $\psi(\ccz{Q}) = (D_\ell/H_\ell) \ee \pm O(\ee^3) = \ee/(\ell+1) \pm O(\ee^3)$.
  Hence, $\Edist{\ccz{Q}}{\ccz{S}} = \ee/(\ell+1) \pm O(\ee^3)$, which implies
  \begin{align}
    D_{Q,S}^2  &=  \frac{1}{(\ell+1)^2} \ee^2 \pm O(\ee^3) ; \\
    R_S^2  &=  R_Q^2 + \frac{1}{(\ell+1)^2} \ee^2 \pm O(\ee^3)
            =  R_\ell^2 + \frac{j+1}{(\ell+1)^2} \ee^2 \pm O(\ee^3),
  \end{align}
  where, to obtain the bounds for $R_S^2$, we used the inductive assumption for $R_Q^2$. 
  This proves \eqref{eqn:L45A} and \eqref{eqn:L45B}.
  To verify \eqref{eqn:L45C}, we note that $X_Q = O(\Delta^3)$ by Lemma~\ref{lem:pyramid_step}.
  The set of points in $\aff{S}$ whose projections to $\aff{Q}$ are at distance at most $X_Q$ from $\ccz{Q}$ is a solid cylinder whose central axis is a line normal to $\aff{Q}$.
  The edge with endpoints $a$ and $a''$ is almost parallel to this axis, so the bisector of the two points intersects the axis almost orthogonally, and certainly at an angle larger than $\pi/4$.
  The points at distance at most $X_S$ from $\ccz{S}$ are contained in the intersection of the cylinder with the slab of points at distance at most $X_E$ from the bisector, which is contained in a ball of radius $\sqrt{2} X_Q + X_E = O(\Delta^3)$.
\end{proof}

This completes the inductive argument, establishing Hypotheses~I, II, and III:  the base case is covered by Lemmas~\ref{lem:bounds_for_long_edges_in_odd_dimensions} and \ref{lem:bounds_for_bisectors_in_odd_dimensions}, and the remaining cases are reached via the two kinds of inductive steps proved in Lemmas~\ref{lem:pyramid_step} and \ref{lem:bi-pyramid_step}.
In particular, the bounds furnished for the $D_{Q,S}$ imply the required bound for $D_S$, which is the minimum over all facets $Q$ of $S$.

\subsection{All Simplices are Critical}
\label{sec:4.3}

The above analysis implies that for sufficiently small $\Delta > 0$ the circumcenter of every simplex in $\Delaunay{A}$ is contained in the interior of the simplex.
This is half of the proof that all simplices in $\Delaunay{A}$ are critical.
The second half of the proof is not difficult.
\begin{corollary}[All Critical in $\Rspace^{2k+1}$]
  \label{cor:all_critical_in_odd_dimensions}
  Let $d = 2k+1$, $n \geq 2$, $\Delta > 0$ sufficiently small, and $A = A_d (n,\Delta) \subseteq \Rspace^d$.
  Then every simplex in $\Delaunay{A}$ is a critical simplex of $\RadiusF \colon \Delaunay{A} \to \Rspace$.
\end{corollary}
\begin{proof}
  A simplex $S \in \Delaunay{A}$ is a critical simplex of $\RadiusF$ iff it contains the circumcenter in its interior, and the $(d-1)$-sphere centered at the circumcenter and passing through the vertices of $S$ does not enclose or pass through any of the other points of $A$.
  By Hypotheses~II and III, the first condition holds.
  To derive a contradiction, assume the second condition fails for $S \in \Delaunay{A}$.
  In other words, there is a point, $b \in A$, that is not a vertex of $S$ but it is enclosed by or lies on the said $(d-1)$-sphere.
  Then $\dime{S} < d$, else the $(d-1)$-sphere intersects each circle in two points, so there is no possibility for another point to interfere.

  Since the $(d-1)$-sphere intersects every circle in only two points, we may assume that $b$ lies on a circle not touched by $S$, or that $b$ neighbors a vertex of $S$ along their circle, and this is the only vertex of $S$ on this circle.
  Then we can add $b$ as a new vertex to get a simplex $T$ with $\dime{T} = \dime{S} + 1$.
  This simplex also belongs to $\Delaunay{A}$, but its circumcenter does not lie in its interior, which contradicts Hypotheses~II and III.
\end{proof}

\subsection{Counting the Cycles}
\label{sec:4.4}

The final counting argument is similar to the one for even dimensions, with a few crucial differences.
Instead of congruent simplices, we have almost congruent simplices in odd dimensions, but they are similar enough to be separated by their circumradii.
\begin{corollary}[Ordering of Radii in $\Rspace^{2k+1}$]
  \label{cor:ordering_of_radii_in_odd_dimensions}
  Let $d = 2k+1$, $n \geq 2$, $\Delta > 0$ sufficiently small, $A = A_{2k+1}(n,\Delta) \subseteq \Rspace^{2k+1}$, and $\RadiusF \colon \Delaunay{A} \to \Rspace$ the radius function.
  Then the circumradii of two simplices, $S, T \in \Delaunay{A}$, satisfy
  $\RadiusF(S) < \RadiusF(T)$ if $\touch{S} < \touch{T}$, or $\touch{S} = \touch{T}$ and $\short{S} < \short{T}$.
\end{corollary}
\begin{proof}
  By Corollary~\ref{cor:all_critical_in_odd_dimensions}, the circumradii are the values of the simplices under the radius function,
  and by Hypothesis~I, the circumradii are segregated into groups according to the number of touched circles and the number of short edges.
  It follows that the values of $\RadiusF$ are segregated the same way.
\end{proof}

We are interested in three kinds of thresholds: the $\varrho_{\ell-1,\ell-1}$, which separate the simplices that touch at most $\ell$ circles from those that touch at least $\ell+1$ circles, the $\varrho_{\ell,-1}$, which separate the $\ell$-simplices without short edges from the other simplices that touch the same number of circles, and the $\varrho_{k,j}$, which separate the $(k+j+1)$-simplices that touch all $k+1$ circles from the $(k+j+2)$-simplices that touch all $k+1$ circles.
We first study the Alpha complexes defined by the first type of thresholds, $\Alfa{\ell-1}{\ell-1} = \RadiusF^{-1} [0, \varrho_{\ell-1,\ell-1}]$.
\begin{lemma}[Constant Homology in $\Rspace^{2k+1}$]
  \label{lem:constant_homology_in_odd_dimensions}
  Let $d = 2k+1$ be a constant, $A = A_{d} (n, \Delta) \subseteq \Rspace^{2k+1}$, and $1 \leq \ell \leq k$.
  Then $\Betti{p}{(\Alfa{\ell-1}{\ell-1})} = O(1)$ for every $p$.
\end{lemma}
\begin{proof}
  Pick $\ell$ of the $k+1$ circles used in the construction of $A$, let $A' \subseteq A$ be the points on these $\ell$ circles, and note that the full subcomplex of $\Delaunay{A}$ with vertices in $A'$ has no non-trivial (reduced) homology.
  We may collapse this full subcomplex to a single $(\ell-1)$-simplex, e.g.\ the $(\ell-1)$-dimensional face of $\Simplex{}$ whose vertices correspond to the $\ell$ circles.

  $\Alfa{\ell-1}{\ell-1}$ is the union of $\binom{k+1}{\ell}$ such full subcomplexes of $\Delaunay{A}$, one for each choice of $\ell$ circles.
  The intersections of these subcomplexes are of the same type, namely induced subcomplexes of $\Delaunay{A}$ for points on $\ell$ or fewer of the circles.
  Hence, $\Alfa{\ell-1}{\ell-1}$ has the homotopy type of the complete $(\ell-1)$-dimensional simplicial complex with $k+1$ vertices, which has a single non-trivial homology group of rank is $\binom{k}{\ell}$.
  As required, this rank is a constant independent of $n$ and $\Delta$.
\end{proof}

To prove relation \eqref{eqn:thm41A} of Theorem~\ref{thm:maximum_Betti_numbers_in_odd_dimensions}, we second consider the Alpha complexes defined by the second type of thresholds, $\Alfa{\ell}{-1} = \RadiusF^{-1} [0, \varrho_{\ell,-1}]$.
This complex is $\Alfa{\ell-1}{\ell-1}$ together with all $\ell$-simplices without short edges.
By Lemma~\ref{lem:constant_homology_in_odd_dimensions}, only a constant number of them give death to $(\ell-1)$-cycles, while all others give birth to $\ell$-cycles.
This implies that the rank of the $\ell$-th homology group of $\Alfa{\ell}{-1}$ is the number of $\ell$-simplices without short edges minus a constant, which is $\binom{k+1}{\ell+1} (n+1)^{\ell+1} \pm O(1)$.
This construction works for $0 \leq \ell \leq k$, which implies relation \eqref{eqn:thm41A}.

\smallskip
To prove relation \eqref{eqn:thm41B} inductively, we third consider the Alpha complexes defined by the third type of thresholds, $\Alfa{k}{j} = \RadiusF^{-1} [0, \varrho_{k,j}]$, for $0 \leq j \leq k$.
The induction hypothesis is
\begin{align}
  \Betti{p}{(\Alfa{k}{p-k-1})}  &=  \tbinom{k}{p-k} \cdot (n+1)^{k+1} \pm O(n^k) ,
\end{align}
and we use the case $p=k$ of relation \eqref{eqn:thm41A} as the induction basis.
The difference between $\Alfa{k}{p-k-1}$ and $\Alfa{k}{p-k}$ are the $(p+1)$-simplices with $p-k+1$ short edges.
Their number is 
\begin{align}
  \tbinom{k+1}{p-k+1} \cdot (n+1)^{2k-p} n^{p-k+1}
    &=  \tbinom{k+1}{p-k+1} \cdot (n+1)^{k+1} \pm O(n^k) ,
\end{align}
This number divides up into the ones that give death and the remaining ones that give birth.
Since $\binom{k+1}{p-k+1} - \binom{k}{p-k} = \binom{k}{p-k+1}$, this implies
\begin{align}
  \Betti{p+1}{(\Alfa{k}{p-k})}  &=  \tbinom{k}{p-k+1} \cdot (n+1)^{k+1} \pm O(n^k) ,
\end{align}
as needed to finish the inductive argument.

\subsection{Voids in Even Dimensions}
\label{sec:4.5}

We return to the one case in $d = 2k$ dimensions that is not covered by the construction in Section~\ref{sec:2}, namely the $(2k-1)$-st Betti number.
It counts the top-dimensional holes, which we refer to as \emph{voids}.
Notwithstanding that the construction in Section~\ref{sec:2} does not provide any voids, Theorem~\ref{thm:maximum_Betti_numbers_in_even_dimensions} claims the existence of $N = k(n+1) + 2$ points in $\Rspace^{2k}$ and a radius such that $\Betti{2k-1}{} = n^k \pm O(n^{k-1})$.

\smallskip
The set of $N$ points whose \v{C}ech complex has that many voids is a straightforward modification of the construction in $2k-1$ dimensions:  place $A = A_{2k-1} (n, \Delta)$ in the $(2k-1)$-dimensional hyperplane $x_{2k} = 0$ in $\Rspace^{2k}$.
Every $(2k-2)$-cycle---which corresponds to a void in $2k-1$ dimensions---is now a pore in the hyperplane that connects the two half-spaces.
In the odd-dimensional construction, all pores arise when the radius is roughly $R_{k-1} \geq \frac{1}{2}$, and they are located in a small neighborhood of the origin.
By choosing $\Delta > 0$ sufficiently small, we can make this neighborhood arbitrarily small.
It is thus easy to add two points, one on each side of the hyperplane, such that their balls close the pores from both sides and turn them into voids in $\Rspace^{2k}$.
More formally, the two points doubly suspend each $(2k-2)$-cycle into a $(2k-1)$-cycle.
Hence, Theorem~\ref{thm:maximum_Betti_numbers_in_odd_dimensions} for $d=2k-1$ and $p=2k-2$, which gives $\Betti{p}{} = (n+1)^k \pm O(n^{k-1})$, provides the missing case in the proof of Theorem~\ref{thm:maximum_Betti_numbers_in_even_dimensions}.

\section{Discussion}
\label{sec:5}

In this paper, we give asymptotically tight bounds for the maximum $p$-th Betti number of the \v{C}ech complex of $N$ points in $\Rspace^d$.
These bounds also apply to the related Alpha complex and the dual union of equal-size balls in $\Rspace^d$.
They do not apply to the Vietoris--Rips complex, which is the flag complex that shares the $1$-skeleton with the \v{C}ech complex for the same data.
In other words, the Vietoris--Rips complex can be constructed by adding all $2$- and higher-dimensional simplices whose complete set of edges belongs the $1$-skeleton of the \v{C}ech complex.
This implies $\Betti{1} (\VRips{A}{r}) \leq \Betti{1} (\Cech{A}{r})$, since adding a triangle may lower but cannot increase the first Betti number.

\smallskip
As proved by Goff~\cite{Gof11}, the 1-st Betti number of the Vietoris--Rips complex of $N$ points is $O(N)$, for all radii and in all dimensions, so also in $\Rspace^3$.
Compare this with the quadratic lower bound for \v{C}ech complexes proved in this paper.
This implies that the first homology group of this \v{C}ech complex has a basis in which most generators are tri-gons; that is: the three edges of a triangle.
The circumradius of a tri-gon is less than $\sqrt{2}$ times the half-length of its longest edge, which implies that most of the $\Theta (N^2)$ generators exist only for a short range of radii.
In the language of persistent homology~\cite{EdHa10}, most points in the $1$-dimensional persistence diagram represent $1$-cycles with small persistence.
Similarly, the $2$-nd Betti number of a Vietoris--Rips complex in $\Rspace^3$ is $o(N^2)$~\cite{Gof11}, compared to that of a \v{C}ech complex, which can be $\Theta (N^2)$.
Hence, most points in the corresponding persistence diagram represent $2$-cycles with small persistence.

\subsection*{Acknowledgments}

\footnotesize{
  The authors thanks Matt Kahle for communicating the question about extremal \v{C}ech complexes, Ben Schweinhart for early discussions on the linked circles construction in three dimensions, and G\'abor Tardos for helpful remarks and suggestions.
}


\appendix
\small

\end{document}